\documentclass[11pt]{article}
\input amssymb.sty
\input amssym.def
\input amssym
\input epsf

\newtheorem{theorem}{Theorem}[section]
\newtheorem{lemma}[theorem]{Lemma}
\newtheorem{proposition}[theorem]{Proposition}

\newtheorem{conjecture}[theorem]{Conjecture}
\newtheorem{definition}[theorem]{Definition}

\newtheorem{theorem-construction}[theorem]{Theorem - Construction}
\newtheorem{lemma-definition}[theorem]{Lemma - Definition}

\begin{document}

\newcommand{\Z}{{\Bbb Z}} 
\newcommand{\R}{{\Bbb R}}
\newcommand{\Q}{{\Bbb Q}}
\newcommand{\C}{{\Bbb C}}
\newcommand{\lra}{\longrightarrow}
\newcommand{\lms}{\longmapsto}
\newcommand{\AAA}{{\Bbb A}}
\newcommand{\Alt}{{\rm Alt}}
\newcommand{\wg}{\wedge}
\newcommand{\ol}{\overline}
\newcommand{\CP}{{\Bbb C}P}
\newcommand{\bwg}{\bigwedge}
\newcommand{\caL}{{\cal L}}
\newcommand{\PP}{{\Bbb P}}
\newcommand{\HH}{{\Bbb H}}
\newcommand{\LL}{{\Bbb L}}

\begin{titlepage}
\title{Regulators}
\author{A. B.  Goncharov}
\date{}
\end{titlepage}
\maketitle

\qquad \qquad \qquad \qquad \qquad \quad {\it To Steve Lichtenbaum  for his 65th birthday}
\tableofcontents

\section{Introduction}

The $\zeta$--function is one of   the most deep and mysterious 
objects in mathematics. During the last two centuries 
it has served as a key source of new ideas and concepts in 
arithmetic algebraic geometry. The 
$\zeta$--function seems to be created to guide  mathematicians 
into the right directions. 
To illustrate this, let me recall three themes  in the 20th century  
mathematics which emerged from the study of the most basic properties of 
$\zeta$--functions: their zeros, analytic properties and special values.

1. Weil's conjectures on $\zeta$--functions 
of varieties over finite fields  inspired  
Grothendieck's revolution in  algebraic geometry 
and led Grothendieck to the concept of motives, and Deligne to 
the  yoga  of weight filtrations. In fact  (pure) motives over $\Q$ can be viewed as
the simplest pieces of algebraic varieties for which the $L$-function can be
defined. Conjecturally the $L$-function characterizes a motive. 

2. Langlands' conjectures predict that $n$--dimensional representations of 
the 
Galois group ${\rm Gal}(\overline \Q/\Q)$ 
correspond to automorphic representations 
of $GL(n)_{/\Q}$. 
The relationship between these seemingly unrelated objects 
was manifested by $L$-functions: 
the Artin $L$-function of the Galois representation 
coincides with  the automorphic $L$--function of the 
corresponding representation of $GL(n)$. 

3. Investigation of the behavior of L-functions of arithmetic schemes 
at integer points, culminated in Beilinson's conjectures, led to the discovery of
 the key principles of the theory of mixed motives.

In this survey we elaborate on a single aspect 
of the third theme: regulators. We focus 
on the analytic and geometric aspects
of the story, and explore several different approaches to motivic complexes and 
regulator maps. We neither touch the seminal 
Birch - Swinnerton-Dyer conjecture and the progress made in its direction nor do we 
consider the vast generalization of this conjecture, due Bloch and Kato \cite{BK}.

{\bf 1. Special values of the Riemann $\zeta$-function and their motivic nature}. 
 Euler proved the famous formula for the special values of the Riemann $\zeta$-function 
at positive even integers:
$$
\zeta(2n) = 
(-1)^k\pi^{2k}\frac{2^{2k-1}}{(2k-1)!}\Bigl(\frac{-B_{2k}}{2k}\Bigr)
$$
All attempts to find a similar formula expressing $\zeta(3), \zeta(5), ... $
 via some known quantities failed. The reason became clear only 
in the recent time: First, the special values $\zeta(n)$ are 
periods  of certain elements 
\begin{equation} \label{2.11.04.1}
\zeta^{\cal M}(n) \in {\rm Ext}_{{\cal M}_T(\Z)}^1(\Z(0), \Z(n)), 
\quad n = 2,3,4, ...
\end{equation}
where on the right stays the extension group in the abelian category 
${\cal M}_T(\Z)$ of mixed Tate motives over ${\rm Spec}(\Z)$ (as 
defined in  \cite{DG}). 
Second, the fact that $\zeta^{\cal M}(2n-1)$ are non torsion  elements should 
imply, according to a version of Grothendieck's conjecture on periods, that $\pi$ 
and $\zeta(3), \zeta(5), ... $ are algebraically independent over $\Q$. 
 The analytic manifestation of the motivic nature of the special values 
is the formula
\begin{equation} \label{12q}
\zeta(n) = \int_{0 < t_1 < ... < t_n < 1} \frac{dt_1}{1-t_1} \wedge \frac{dt_2}{t_2} 
\wedge ... \wedge \frac{dt_n}{t_n}
\end{equation}
 discovered 
by  Leibniz. This formula presents $\zeta(n)$ as a length $n$ iterated integral. 
The existence of such a formula seems to be a specific property 
of the  $L$-values at integer points. 
A geometric construction of the motivic $\zeta$-element (\ref{2.11.04.1}) 
using the moduli space ${\cal M}_{0, n+3}$ 
is given in Chapter 4.5.

Beilinson \cite{B1} conjectured a similar 
picture for special values of $L$-functions of motives at integer points. 
In particular, his conjectures imply that these special values should be
periods  
(in fact a very special kind of periods). 
For the Riemann $\zeta$-function this is given by the formula (\ref{12q}). 
In general a period is a number given by an integral
$$
\int_{\Delta_B}\Omega_A
$$ 
where $\Omega_A$ is a differential form on a variety $X$ with singularities at 
a divisor $A$, $\Delta_B$ is a chain with boundary at a divisor $B$, and   
$X$, $A$, $B$ are defined over $\Q$. So far we can write    
$L$-values at integer points 
as  periods only in a few cases. Nevertheless, in all cases when Beilinson's conjecture 
was confirmed,  we have such 
a presentation.  More specifically, such a presentation for the special values 
of the Dedekind $\zeta$-function this comes from the Tamagawa Number formula and
Borel's work \cite{Bo2}, and in the other cases it is given by
Rankin-Selberg type formulas. 
In general the mechanism 
staying behind this phenomenon remains a mystery.

Let us now turn to another classical example: the residue of the Dedekind 
$\zeta$-function at $s=1$.


{\bf 2. The class number formula and 
the weight one Arakelov motivic complex}.  
Let $F$ be a number field with $r_1$ real and $r_2$ complex places, so that 
$[F:\Q] = 2r_1+r_2$. Let $\zeta_F(s)$ be the Dedekind 
$\zeta$-function of $F$. Then according to  Dirichlet and  Dedekind one has 
\begin{equation} \label{3.20.03.1}
{\rm Res}_{s=1}\zeta_F(s) = \frac{2^{r_1+r_2}\pi^{r_2}R_F h_F}{w_F |D_F|^{1/2}}
\end{equation}
Here $w_F$ is the number of roots 
of unity in $F$, $D_F$ is the discriminant, 
$h_F$ is the class number, 
and $R_F$ is the regulator of $F$, whose definition we recall below. 
Using the functional equation 
for $\zeta_F(s)$, S. Lichtenbaum \cite{Li1} wrote (\ref{3.20.03.1}) as 
\begin{equation} \label{3.21.03.1}
\lim_{s \to 0}s^{-(r_1+r_2-1)}\zeta_F(s) = -\frac{h_FR_F}{w_F}
\end{equation}
Let us interpret the right hand side of this formula 
 via  the weight one Arakelov 
motivic complex for ${\rm Spec}({\cal O}_F)$, where 
${\cal O}_F$ is the ring of integers in $F$. 
Let us define  the following diagram, where ${\cal P}$ runs 
through all prime ideals of ${\cal O}_F$:
\begin{equation} \label{3.21.03.5}
\begin{array}{ccc}
\R^{r_1+r_2}& \stackrel{\Sigma}{\lra} & \R\\
&&\\
R_1\uparrow &&\uparrow l\\
&&\\
F^*&\stackrel{\rm div}{\lra}&\oplus_{\cal P}\Z
\end{array}
\end{equation}
In this diagram, the maps are given as follows: If ${\rm val}_{{\cal P}}$ 
is the canonical valuation 
defined on $F$ by  ${\cal P}$ and  $|{\cal P}|$ is 
the norm of ${\cal P}$, then 
$$
{\rm div}(x) = \sum {\rm val}_{{\cal P}}(x)[{\cal P}], 
\qquad l: [{\cal P}] \lms -\log |{\cal P}|, \qquad \Sigma: (x_1, ..., x_n) 
\lms  
\Sigma x_i
$$ 
The regulator map $R_1$ is defined by 
$
x \in F^* \lms (\log |x|_{\sigma_1}, ..., \log |x|_{\sigma_{r_1+r_2}})
$, 
where $\{\sigma_1, ..., \sigma_{r_1+r_2}\}$ 
is the set of all archimedian places of $F$ and 
$|\ast|_{\sigma}$ is the valuation defined 
$\sigma$, 
($|x|_{\sigma}:= |\sigma(x)|^2$ for a complex place $\sigma$).

The product formula tells  us $\Sigma \circ R_1 + l \circ {\rm div} = 0$. 
Therefore summing up the groups over the diagonals in (\ref{3.21.03.5})  
we get a complex. 
The first two groups of this complex,  placed in degrees $[1,2]$, 
form   the {\it weight 
one Arakelov motivic complex} ${\Gamma}_{\cal A}({\cal O}_F, 1)$ 
of ${\cal O}_F$. There is  a map 
$$
H^2{\Gamma}_{\cal A}({\cal O}_F, 1) \stackrel{}{\lra} \R
$$ 
Let $\widetilde H^2{\Gamma}_{\cal A}({\cal O}_F, 1)$ be its kernel. 
Then there is an exact sequence 
$$
0 \lra \frac{{\rm Ker}(\Sigma)}{R_1({\cal O}^*_F)}\lra \widetilde H^2 {\Gamma}_{\cal A}({\cal O}_F, 1) \lra {\rm Cl}_F \lra 0 
$$
Further, $R_1({\cal O}^*_F)$ is a lattice in 
${\rm Ker}(\Sigma) = \R^{r_1+r_2-1}$, its volume with respect to 
the measure $\delta(\sum x_i) dx_1\wedge ... \wedge dx_{r_1+r_2}$ 
is $R_F$, and $h_F = |{\rm Cl}_F|$. Therefore 
$$
{\rm vol}\widetilde H^2 {\Gamma}_{\cal A}({\cal O}_F, 1) =  R_F h_F; \quad 
H^1\widetilde {\Gamma}_{\cal A}({\cal O}_F, 1) = \mu_F
$$
Now the class number formula (\ref{3.21.03.1}) reads
\begin{equation} \label{3.03.21.6}
\lim_{s \to 0}s^{-(r_1+r_2-1)}\zeta_F(s) = -
\frac{{\rm vol}\widetilde H^2 {\Gamma}_{\cal A}({\cal O}_F, 1)}{H^1{\Gamma}_{\cal A}({\cal O}_F, 1) } 
\end{equation}
The right hand side is a volume of the determinant of a complex, see 
 Chapter 2.5. I do not know the cohomological 
origin of the sign  in (\ref{3.03.21.6}).

 {\bf 3. Special values of the Dedekind $\zeta$-functions, 
Borel regulators and polylogarithms}. B. Birch and J. 
Tate \cite{T} proposed a generalization of the class number formula 
for totally real fields using Milnor's $K_2$--group 
of  ${\cal O}_F$: 
$$
\zeta_F(-1) = \pm \frac{|K_2({\cal O}_F)|}{w_2(F)}
$$
Here $w_2(F)$ is the largest integer  $m$ such that  
${\rm Gal}(\overline F/F)$ acts trivially on $\mu_{\infty}^{\otimes 2}$.  
Up to a power of $2$, the above formula follows from the Iwasawa main conjecture 
for totally real fields, 
proved by B. Mazur and A. Wiles for $\Q$ \cite{MW}, and by A. Wiles \cite{W}
 in general.

S. Lichtenbaum \cite{Li1} suggested that for 
$\zeta_F(n)$ there should be a  formula similar to (\ref{3.21.03.1}) with a  
higher regulator defined using Quillen's $K$-groups 
$K_{*}(F)$ of 
$F$. Such a formula for $\zeta_F(n)$, considered up to a 
non zero rational factor, has been established soon after 
in the fundamental work of A. Borel \cite{Bo1}-\cite{Bo2}. 
Let us discuss it in more detail. The rational $K$--groups of a field 
can be defined as the primitive part in the homology of $GL$. 
Even better,  one can show that 
$$
K_{2n-1}(F) \otimes \Q \stackrel{\sim}{=} {\rm Prim}H_{2n-1}(GL_{2n-1}(F), \Q)
$$
Let $\R(n):= (2\pi i)^n\R$. There is a distinguished 
class, called the Borel class, 
$$
B_n \in H^{2n-1}_c(GL_{2n-1}(\C), \R(n-1))
$$
in the continuous cohomology of the Lie group $GL_{2n-1}(\C)$. 
Pairing with this class provides the  Borel 
regulator map 
$$
r^{\rm Bo}_n: K_{2n-1}(\C) \lra \R(n-1)
$$
Let $X_F := \Z^{Hom(F, \C)}$. 
The Borel regulator map on $K_{2n-1}(F)$ is the composition 
$$
K_{2n-1}(F) \lra \oplus_{Hom(F, \C)} K_{2n-1}(\C) 
 \lra X_F \otimes \R(n-1)
$$
 The image of this map 
is invariant under complex conjugation acting  both  on 
$Hom(F, \C)$ and $\R(n-1)$. So we get the map 
\begin{equation} \label{3.19.03.2}
R^{\rm Bo}_n: K_{2n-1}(F) \lra \left(X_F \otimes \R(n-1)  \right)^+
\end{equation}
Here $+$ means the invariants under the complex conjugation. 
Borel proved that for $n>1$ the image of this map is a lattice, and 
the volume $R_n(F)$ of this lattice is related to 
the Dedekind $\zeta$-function as follows:
$$
R_n(F) \sim_{\Q^*}\lim_{s \to 1-n}(s-1+n)^{-d_n} \zeta_F(s), \qquad 
$$
Here $a \sim_{\Q^*}b$ means $a = \lambda b$ for some $\lambda \in \Q^*$, and 
$$
d_n = {\rm dim} \left(X_F \otimes \R(n-1)  \right)^+
= \left\{ \begin{array}{ll} 
r_1+r_2&: \mbox{$n>1$ odd} \\ 
r_2 &: \mbox{$n\geq 2$ even} \end{array}\right.
 $$ 
Using the functional equation for $\zeta_F(s)$ it 
tells us about $\sim_{\Q^*}\zeta_F(n)$. 
However Lichtenbaum's original conjecture 
was stronger since it was  about $\zeta_F(n)$ itself.

In 1977 S. Bloch discovered \cite{Bl4}, \cite{Bl5} that 
the regulator map on $K_3(\C)$ can be explicitly defined 
using the dilogarithm. Here is how the story looks today. 
The dilogarithm \index{dilogarithm} is a multivalued analytic function on 
$\C\PP^1 - \{0, 1, \infty\}$:
$$
{\rm Li}_2(z) := -\int_0^z\log(1-z)\frac{dz}{z}; \qquad 
{\rm Li}_2(z) = \sum_{k=1}^{\infty}\frac{z^2}{k^2} \quad \mbox{for $|z|\leq 1$} 
$$
The dilogarithm has a single-valued version, called the Bloch-Wigner function: 
\index{Bloch-Wigner function}
$$
{\cal L}_2(z) := {\rm Im}\left({\rm Li}_2(z) 
+ \log(1-z)\log|z|\right)
$$
It vanishes on the real line. 
Denote by $r(z_1, ..., z_4)$ the cross-ratio 
of the four points  $z_1, ..., z_4$ on the projective line. 
The Bloch-Wigner function satisfies Abel's five term relation: 
\index{Abel's five term relation}
for any five points $z_1, ... ..., z_5$ on $\C\PP^1$ one has 
\begin{equation} \label{3.28.03.1}
\sum_{i=1}^5(-1)^i{\cal L}_2(r(z_1, ..., \widehat z_i, ..., z_5)) = 0
\end{equation}

Let ${\cal H}_3$ be the hyperbolic three space. Its absolute is identified 
with the Riemann sphere $\C\PP^1$. Let $I(z_1, ..., z_4)$ 
be the ideal geodesic simplex 
with the vertices at the points $z_1, ..., z_4$ at the absolute. 
Lobachevsky proved that 
$$
{\rm vol}I(z_1, ..., z_4) = {\cal L}_2(r(z_1, ..., z_4))
$$
(He got this in a different but equivalent form). 
Lobachevsky's formula makes Abel's equation obvious: the alternating sum 
of the geodesic simplices with the 
vertices at  $z_1, ..., \widehat z_i, ..., z_5$ is empty.

Abel's equation  can be interpreted 
as follows: for any $z \in \C\PP^1$ the function
\begin{equation} \label{3.17.03.2}
{\cal L}_2(r(g_1 z, ..., g_4 z)), \quad g_i \in GL_2(\C)
\end{equation}
is a measurable $3$-cocycle of the Lie group $GL_2(\C)$. 
The cohomology class of this cocycle is non trivial. 
The simplest way to see it is this. 
The function 
$$
{\rm vol}I(g_1x, ..., g_4x), \quad g_i \in GL_2(\C), x \in {\cal H}^3
$$
provides a smooth $3$-cocycle of  $GL_2(\C)$. Its  cohomology class is 
nontrivial: indeed,  its infinitesimal version is provided by the volume form 
in ${\cal H}_3$. Since 
the cohomology class does not depend on $x$, we can 
take $x$ at the absolute, proving the claim.   


Let $F$ be an arbitrary field. Denote by 
$\Z[F^*]$ the free abelian group generated 
by the set $F^*$. Let $R_2(F)$ be the subgroup of $\Z[F^*]$ 
generated by the elements
$$
\sum_{i=1}^5(-1)^i\{r(z_1, ..., \widehat z_i, ..., z_5)\}
$$
Let $B_2(F)$ be the quotient of $\Z[F^*]$ by the subgroup $R_2(F)$.
Then one shows that the map 
$$
\Z[F^*] \lra \Lambda^2F^*; \quad \{z\} \lms (1-z) \wedge z
$$
kills the subgroup $R_2(F)$, providing a complex 
(called the Bloch-Suslin complex) 
\begin{equation} \label{3.31.03.2}
\delta_2: B_2(F) \lra \Lambda^2F^*
\end{equation}
By Matsumoto's theorem ${\rm Coker} \delta_2 = K_2(F)$. 
Let us define the Milnor ring $K^M_*(F)$ of $F$ as the quotient of the
tensor algebra of the abelian group $F^*$ by the two sided ideal generated by
the Steinberg elements $(1-x) \otimes x$ where $x \in F^*-1$. The product map in the $K$-theory 
provides a map 
$$
\otimes^nK_1(F) = \otimes^nF^* \lra K_n(F)
$$
It kills the Steinberg elements, and thus provides a map $K_n^M(F) \to
K_n(F)$. 
Set 
$$
K_3^{\rm ind}(F):= {\rm Coker}(K_3^M(F) \lra K_3(F))
$$ 
A.A. Suslin \cite{Su} proved that there is an exact sequence 
\begin{equation} \label{3.17.03.1}
0 \lra {\rm Tor}(F^*, F^*)^{\sim} \lra K_3^{\rm ind}(F) \lra {\rm Ker}\delta_2
\lra 0
\end{equation}
where ${\rm Tor}(F^*, F^*)^{\sim}$ is a nontrivial extension 
of $\Z/2\Z$ by ${\rm Tor}(F^*, F^*)$. 

Abel's relation provides  a well defined 
homomorphism
$$
{\cal L}_2: B_2(\C) \lra \R; \quad \{z\}_2 \lms {\cal L}_2(z)
$$
Restricting it 
to the subgroup ${\rm Ker}\delta_2 \subset B_2(\C)$,  
and using  (\ref{3.17.03.1}),  
we get a map  
$K_3^{\rm ind}(\C) \lra \R$. Using the interpretation 
of the cohomology class of the cocycle (\ref{3.17.03.2}) 
as a volume of geodesic simplex,  
one can show that it is essentially the 
Borel regulator. Combining this 
with Borel's theorem we get an explicit formula for 
$\zeta_F(2)$ for an arbitrary number field $F$. 

How to  generalize this beautiful story? Recall the classical polylogarithms 
\index{classical polylogarithms}
$$
{\rm Li}_n(z) = \sum_{k=1}^{\infty}\frac{z^k}{k^n}, \quad |z| \leq 1; \qquad 
{\rm Li}_n(z) = \int_0^z{\rm Li}_{n-1}(z)\frac{dz}{z}
$$
D. Zagier \cite{Z1} formulated a precise conjecture 
expressing $\zeta_F(n)$ via classical polylogarithms, 
see the survey \cite{GaZ}. 
It was proved for $n=3$ in \cite{G1}-\cite{G2}, but  
it is not known for $n>3$, although its easier part 
has been proved in \cite{dJ3}, \cite{BD2} and, in a different way, 
 in Chapter 4.4 below.

Most of the $\zeta_F(2)$ picture  has been 
generalized to the case of $\zeta_F(3)$ in \cite{G1}-\cite{G2} and \cite{G3}. 
Namely, there is a single valued version of the trilogarithm: 
$$
{\cal L}_3(z):= {\rm Re}\left({\rm Li}_3(z) - {\rm Li}_2(z)\log|z| 
+ \frac{1}{6}{\rm Li}_1(z)\log^2|z| \right)
$$
It satisfies the following functional equation which generalizes
(\ref{3.28.03.1}). Let us define the generalized cross-ratio 
of 6 points $x_0, ..., x_5$ in $P^2$ \index{generalized cross-ratio 
of 6 points in $P^2$} as follows.  
We present $P^2$ as a projectivization of the three dimensional vector 
space $V_3$ and choose the vectors $l_i \in V_3$ projecting to the points $x_i$. 
Let us choose a volume form $\omega \in {\rm det}V^*_3$ and set  
$\Delta(a,b,c):= <a \wedge b \wedge c, \omega>$. Set
$$
r_3(x_0, ..., x_5):= {\rm Alt}_6
\{ \frac{\Delta(l_0, l_1, l_3)\Delta(l_1, l_2, l_4)
\Delta( l_2, l_0, l_5)}{\Delta(l_0, l_1, l_4)
\Delta(l_1, l_2, l_5)\Delta(l_2, l_0, l_3)} \} \in \Z[F^*]
$$
Here ${\rm Alt}_6$ denotes the alternation of $l_0, ..., l_5$. 
The function ${\cal L}_3$ extends by linearity to
 a homomorphism ${\cal L}_3: \Z[\C^*] \lra \R$, 
and there is a generalization of Abel's equation to the 
case of the trilogarithm (see \cite{G2} and the appendix to \cite{G3}):
\begin{equation} \label{3.28.03.2}
\sum_{i=1}^7(-1)^i{\cal L}_3(r_3(x_1, ..., \widehat x_i, ... x_7)) = 0
\end{equation}
We define the group $B_3(F)$ as the quotient 
of $\Z[F^*]$ by the subgroup generated by  the functional equations 
(\ref{3.28.03.2})
for the trilogarithm. \index{functional equations 
for the trilogarithm} Then there is  complex 
\begin{equation} \label{3.19.03.1}
B_3(F) \stackrel{\delta_3}{\lra} B_2(F)\otimes F^* 
\stackrel{\delta_2 \wedge {\rm Id}}{\lra} \Lambda^3F^*; 
\quad \delta_3: \{x\}_3 \lms \{x\}_2 \otimes x
\end{equation}
There is a map $K_5(F) \lra {\rm Ker}\delta_3$ such that in the case 
$F = \C$ the composition
$$
K_5(\C) \lra {\rm Ker}\delta_3 \hookrightarrow B_3(\C) 
\stackrel{{\cal L}_3}{\lra} \R
$$
coincides with the Borel regulator (see \cite{G2} and appendix in \cite{G3}). 
This plus Borel's theorem 
leads to an explicit formula expressing  $\zeta_F(3)$ via 
the trilogarithm conjectured by 
Zagier \cite{Z1}.

In Chapter 3 we define, following \cite{G4}, the Grassmannian 
$n$-logarithm function ${\cal L}^G_n$. 
It is a function on the configurations of $2n$ 
hyperplanes in $\C\PP^{n-1}$. 
One of its functional equations generalizes Abel's equation: 
$$
\sum_{i=1}^{2n+1}(-1)^i{\cal L}^G_n(h_1, ..., 
\widehat h_i, ..., h_{2n+1})=0
$$
It means that for a given hyperplane $h$ the function 
${\cal L}^G_n(g_1 h, ...,  g_{2n}h)$, $g_i \in GL_n(\C)$, 
is a  measurable cocycle of $GL_n(\C)$. Its 
cohomology class essentially coincides with the restriction of 
the Borel class $B_n$ to $GL_n(\C)$. 
To prove this we  show that   ${\cal L}^G_n$ 
 is the boundary 
value of a certain function  
defined on configurations of $2n$ 
points in the symmetric space $SL_n(\C)/SU(n)$.  
Using this we express the Borel regulator via the 
Grassmannian polylogarithm. 
We show that for $n=2$ we recover the dilogarithm story: 
 ${\cal L}^G_2$ coincides with the Bloch-Wigner function,  
$SL_2(\C)/SU(2)$ 
is the hyperbolic space, and the extension of ${\cal L}^G_2$ 
is given by the volume of geodesic simplices.  The proofs can be found in 
\cite{G7}.

{\bf 4. Beilinson's conjectures and Arakelov motivic complexes}. 
A conjectural generalization of the class number formula (\ref{3.20.03.1}) 
to the case of elliptic curves was suggested in the seminal work of 
Birch and Swinnerton-Dyer. Several years later J.Tate 
formulated conjectures relating algebraic cycles 
to the poles of zeta functions of algebraic varieties.

Let $X$ be  a regular algebraic variety
over a number field $F$. 
Generalizing the previous works of 
Bloch \cite{Bl4} and P. Deligne \cite{D}, 
A.A. Beilinson \cite{B1} suggested  a fantastic 
picture unifying all the above conjectures.  
 Beilinson defined the rational motivic cohomology of $X$ 
via the algebraic K-theory of $X$ by 
$$
H^i_{Mot}(X, \Q(n)):= gr^{\gamma}_nK_{2n-i}(X) \otimes \Q 
$$
Here $\gamma$ is the Adams $\gamma$-filtration. 

Let us assume that $X$ is projective. 
For schemes which admit 
regular models over $\Z$, Beilinson \cite{B1} defined a 
$\Q$-vector subspace 
\begin{equation} \label{reguula}
H^{i}_{Mot/\Z}(X, \Q(n)) \subset H^{i}_{Mot}(X, \Q(n)),
\end{equation} 
called {\it integral part in the motivic cohomology}. 
Using alterations, A. Scholl 
\cite{Sch} extended this definition to arbitrary regular projective schemes 
over a number field $F$. For every regular, projective and flat over $\Z$ 
model $X'$ of $X$ the subspace (\ref{reguula})  coincides with 
the image of the map $H^{i}_{Mot}(X', \Q(n)) \lra H^{i}_{Mot}(X, \Q(n))$.

For a regular complex projective variety  $X$ 
Beilinson  
constructed in \cite{B1} the regulator map to the Deligne cohomology 
of $X$:
\begin{equation} \label{reguul}
H^i_{Mot}(X, \Q(n)) \lra H^i_{{\cal D}}(X(\C), \Z(n))
\end{equation}
It can be projected to the  Deligne cohomology with real coefficients 
$H^i_{{\cal D}}(X(\C), \R(n))$. 

Now let  $X$ again be a regular projective scheme over a number field $F$, We will view it as a scheme over $\Q$ via the projection 
$X \to {\rm Spec}(F) \to {\rm Spec}(\Q)$. We define the {\it real Deligne cohomology of $X$} 
as the following $\R$-vector space: 
$$
H^{i}_{{\cal D}}
\Bigl((X\otimes_{\Q}\R)(\C), \R(n)\Bigr)^{\overline F_{\infty}}
$$
  where $\overline F_{\infty}$  is an involution given by the composition of 
complex conjugation acting on $X(\C)$ and on the coefficients. Then, 
restricting the map (\ref{reguul}) 
to the integral part in the motivic cohomology (\ref{reguula}) and projecting  
the image onto 
the real Deligne cohomology of $X$, we get a regulator map 
$$
r_{Be}: H^{i}_{Mot/\Z}(X, \Q(n)) \lra H^{i}_{{\cal D}}
\Bigl((X\otimes_{\Q}\R)(\C), \R(n)\Bigr)^{\overline F_{\infty}}
$$
 Beilinson 
formulated a conjecture expressing the special values 
$\sim_{\Q^*}L(h^{i-1}(X), n)$ of the L-function of Grothendieck's motive $h^{i-1}(X)$ 
in terms of this  regulator map, see \cite{B1}, the 
survey \cite{R} and the book \cite{RSS} for the original
 version of the conjecture, and 
the survey by J. Nekov\'a\v{r} \cite {N} for a motivic reformulation.  
 For $X= {\rm Spec}(F)$, 
where $F$ is a number field,  it boils down to 
Borel's theorem. 
A precise Tamagawa Number conjecture about the special values 
$L(h^i(X), n)$ was suggested by Bloch and Kato \cite{BK}.

Beilinson \cite{B2}  and Lichtenbaum \cite{Li2}
 conjectured that the weight $n$ integral 
motivic cohomology of a scheme $X$ should appear 
as cohomology of some complexes of abelian groups ${\Z}_X^{\bullet}(n)$, 
called the weight $n$ motivic complexes of $X$. 
One must have 
$$
H^i_{Mot}(X, \Q(n)) = H^i{\Z}_X^{\bullet}(n)\otimes \Q
$$
Motivic complexes  are objects of the derived category. 
They should appear as hypercohomology of certain complexes of sheaves in Zarisky topology on $X$. 
Beilinson conjectured \cite{B2} that there exists an 
abelian category ${\cal M}{\cal S}_X$ 
of mixed motivic sheaves on $X$,  and that 
one should have an isomorphism in the derived category
\begin{equation} \label{MM}
{\Z}_X^{\bullet}(n) = RHom_{{\cal M}{\cal S}_X}(\Q(0)_X, \Q(n)_X)
 \end{equation}
Here $\Q(n)_X:= p^*\Q(n)$, where $p:X \to {\rm Spec}(F)$ 
is the structure morphism, 
 is a motivic sheaf on $X$ obtained by pull back 
of the Tate motive $\Q(n)$ over the point ${\rm Spec}(F)$. 
This formula implies  the Beilinson-Soul\'e vanishing conjecture: 
$$
H^i_{Mot}(X, \Q(n)) = 0\quad \mbox{for $i<0$ and $i=0, n>0$}
$$ Indeed, 
the negative Ext's between objects of an abelian category are zero, and
we assume that the objects  $\Q(n)_X$ are mutually non-isomorphic. 
Therefore it is quite natural to look for representatives of motivic
complexes which are zero in the negative degrees, as well as in the degree
zero for $n=0$. 

Motivic complexes are more fundamental, and in fact 
simpler  objects then rational $K$--groups. 
Several constructions of motivic complexes are known. 

i) Bloch \cite{Bl1}-\cite{Bl2} suggested a construction 
of the motivic complexes, called the  Higher Chow complexes, 
 using algebraic cycles. The weight $n$ cycle complex 
 appears in a very natural way as a 
``resolution'' for the codimension $n$ Chow groups on $X$ 
modulo rational equivalence, 
see s. 2.1 below.   
These complexes as well as their versions defined by Suslin and Voevodsky 
played an essential role in the construction of 
triangulated categories of mixed motives \cite{V}, \cite{Lev}, \cite{Ha}. 
However  they are unbounded from the left.

ii) Here is a totally different construction of
 the first few of motivic complexes. 
One has ${\Z}_X^{\bullet}(0) := \Z$. Let $X^{(k)}$ 
be the set of all irreducible codimension $k$ subschemes of a scheme $X$.  Then 
$ {\Z}_X^{\bullet}(1)$ is the complex
$$
 {\cal O}^*_X \stackrel{\partial}{\lra} \oplus_{Y \in X^{(1)}}\Z;
\qquad \partial (f):= {\rm div}(f)
$$
The complex $ {\Q}_X^{\bullet}(2)$ is defined as follows. First, 
using the Bloch-Suslin complex (\ref{3.31.03.2}), we define the following complex
$$
\qquad {B}_2(\Q(X)) \stackrel{\delta_2}{\lra} \Lambda^2\Q(X)^*\stackrel{\partial_1}{\lra} \oplus_{Y \in X^{(1)}}\Q(Y)^* \stackrel{\partial_2}{\lra} 
\oplus_{Y \in X^{(2)}}\Z
$$ 
Then tensoring it by $\Q$ we get ${\Q}_X^{\bullet}(2)$. 
Here $\partial_2$ is the tame symbol, and $\partial_1$ 
is given by the divisor of a function on $Y$. 
Similarly one can define a complex $ {\Q}_X^{\bullet}(3)$ using the complex 
(\ref{3.19.03.1}). 
Unlike the cycle complexes, these complexes are
 concentrated exactly in the degrees 
where they  might have nontrivial cohomology. 
It is amazing that  motivic complexes 
have two so different and beautiful incarnations. 
 
Generalizing  this, we introduce in Chapter 4 the {\it polylogarithmic motivic
  complexes}, which are conjectured to be the motivic complexes for an
  arbitrary 
field 
(\cite{G1}-\cite{G2}). Then we  give a motivic proof of 
the weak version of Zagier's conjecture for a number field $F$.  
In Chapter 5 we discuss how to define motivic complexes for an arbitrary
  regular variety $X$ using the polylogarithmic 
motivic complexes of its
  points. 

In Chapters 2 and 5 we discuss constructions of the regulator map  
on the level of complexes. Precisely, we want to define 
for a regular complex projective variety $X$  
a homomorphisms of complexes of abelian groups 
\begin{equation} \label{3/1/00.1}
\mbox{a weight $n$ motivic complex of $X$} \quad \lra \quad 
\mbox{a weight $n$ Deligne complex of $X(\C)$}
\end{equation}
In Chapter 2 we present a construction of a regulator map on the  Higher
Chow complexes 
given in \cite{G4}, \cite{G7}.  
In Chapter 5 we define, following \cite{G6}, a  regulator map 
on polylogarithmic complexes. 
It is given explicitly in terms of the classical polylogarithms.  
Combining this with Beilinson's conjectures we arrive at  
explicit conjectures expressing the  special values of L-functions 
via classical polylogarithms. 
If $X = {\rm Spec}(F)$, where $F$ is a 
number field, this boils down to Zagier's conjecture.

The cone of the map (\ref{3/1/00.1}), shifted by $-1$, 
 defines the {\it weight $n$ Arakelov motivic complex}, and so its cohomology 
 are the {\it weight $n$ Arakelov motivic cohomology} of $X$. 

The weight $n$ Arakelov motivic complex should be considered as an ingredient 
of a definition of the {\it weight $n$ arithmetic motivic complex}. Namely,
one should exist a complex computing the weight $n$ integral motivic cohomology of $X$,
and a natural map from this complex to the weight $n$ Deligne complex. The cone of this
map, shifted by $-1$, would give the weight $n$ arithmetic motivic complex. 
The weight one arithmetic motivic complex  
is the complex $\Gamma_{\cal A}({\cal O}_F, 1)$.

In Chapter 6 we discuss a yet another  approach to motivic 
complexes of fields:    
as standard cochain complexes of the motivic Lie algebras. 
We also discuss a relationship between 
the 
motivic Lie algebra of a field and the (motivic) Grassmannian 
polylogarithms.   

{\bf Acknowledgement}.
This work was partially  supported by 
IHES (Bures sur Yvette), MPI (Bonn), and  the NSF grants 
DMS-0099390. I am grateful to referee for useful comments.

\section  {Arakelov motivic complexes}

In this chapter we define a regulator map from 
the weight $n$ motivic complex, understood as 
Bloch's  Higher Chow groups complex \cite{Bl1}, to 
the weight $n$ Deligne complex.
This map was defined in \cite{G4}, and elaborated in detail in \cite{G7}. 
The construction can be immediately adopted 
to the Suslin-Voevodsky motivic complexes.

{\bf 1.  Bloch's  cycle complex \cite{Bl1}}. \index{Bloch's  cycle complex}
A non degenerate simplex in $\PP^m$ is an ordered
collection of hyperplanes $L_0,...,L_m$  with empty
intersection. 
Let us choose in $\PP^m$ a simplex $L$ and a generic hyperplane $H$.
Then $L$ provides  a simplex in the 
affine space $\AAA^m:= \PP^m  - H$.

Let $X$ be a regular scheme over a field. 
Let $I=(i_1,...,i_k)$ and $L_I:= L_{i_1} \cap ... 
\cap L_{i_k}$. Let ${\cal Z}_{m}(X; n)$ be
the free abelian group generated by irreducible codimension $n$
algebraic 
subvarieties in $X \times \AAA^m$ which  intersect properly (i.e.  the intersection has the  right dimension) all faces
$X \times L_I$. 
Intersection with the codimension 1 face $X \times L_i$ provides a group 
homomorphism  $\partial_i: {\cal Z}_{m}(X;n) \longrightarrow 
{\cal Z}_{m-1}(X; n)$. Set $\partial:= \sum_{i=0}^m (-1)^i \partial_i$. Then   $\partial^2 =0$, so 
$({\cal Z}_{\bullet}(X; n), \partial)$    is a homological complex. 
Consider the  cohomological complex 
${\cal Z}^{\bullet}(X; n):={\cal Z}_{2r-\bullet}(X; n)$. 
Its cohomology give the motivic cohomology of $X$:
$$
H_{{\cal M}}^i(X, \Z(n)):= H^i({\cal Z}^{\bullet}(X; n))
$$
According to the fundamental theorem of Bloch (\cite{Bl1}, \cite{Bl2})
$$
H^i({\cal Z}^{\bullet}(X; n)\otimes \Q)
 = gr_n^{\gamma}K_{2n-i}(X)\otimes \Q
$$

{\bf 2. The Deligne cohomology and Deligne's complex}. \index{Deligne's complex} Let $X$ be a regular
projective variety over $\C$.  
The Beilinson-Deligne complex ${\R}^{\bullet}(X; n)_{{\cal D}}$ 
is the following complex of sheaves in the classical topology on $X(\C)$: 
$$
\underline\R(n) \lra {\cal O}_X \lra \Omega^{1}_{X} 
\lra \Omega^{2}_{X}\lra ... \lra \Omega^{n-1}_{X}
$$
Here the constant sheaf $\underline\R(n) := (2\pi i)^n\underline\R$ is in the degree zero. The
hypercohomology of this complex of sheaves is called the {\it weight $n$ 
Deligne cohomology} \index{weight $n$ 
Deligne cohomology} of
$X(\C)$. They are finite dimensional real vector spaces. 
Beilinson proved \cite{B3} that the truncated 
weight $n$ Deligne cohomology, which are obtained by putting  
the weight $n$  Deligne cohomology equal to zero 
in the degrees $>2n$, can be interpreted as the absolute Hodge cohomology of
$X(\C)$. 

One can replace the above complex of sheaves by a quasiisomorphic one, defined
as the total complex associated with the following
bicomplex: 
$$
\begin{array}{ccccccccccc} \label{del}
\Bigl(\underline{\cal D}_{X}^{0}&\stackrel{d}{\longrightarrow}&\underline{\cal
D}_{X}^{1}&\stackrel{d}{\longrightarrow}&\ldots&\stackrel{d}{\longrightarrow}&\underline{\cal
D}^{n}_{X}&\stackrel{d}{\longrightarrow}&\underline{\cal
D}_{X}^{n+1}&\stackrel{d}{\longrightarrow}&\ldots\Bigr) \otimes \R(n-1)\\
&&&&&&&&&&\\
&&&&&&\uparrow\pi_{n} &&\uparrow\pi_{n}&&\\
&&&&&&&&&&\\
&&&&&&\Omega^{n}_{X}
&\stackrel{\partial}{\longrightarrow}&\Omega_{X}^{n+1}&\stackrel
{\partial}{\longrightarrow}&
\end{array}
$$
Here  $\underline{\cal D}_{X}^{k}$ is the sheaf of real $k$-distributions 
on $X(\C)$, that is $k$-forms with the generalized function coefficients. 
Further, 
$$
\pi_n: \underline{\cal D}_{X}^{p}\otimes \C \longrightarrow
\underline{\cal D}_{X}^{p}\otimes \R(n-1)
$$ is the projection induced by the one  $\C = \R(n-1) \oplus \R(n)  \longrightarrow 
\R(n-1)$, the sheaf  $\underline{\cal D}^{0}_{X}$ is placed in
degree 1, and  
$(\Omega^{\bullet}_{X}, \partial)$ is the De Rham complex of 
sheaves of holomorphic forms. 

To calculate the hypercohomology with coefficients in this complex 
we replace the holomorphic 
de Rham complex by its Doulbeut resolution, take the global
sections of the obtained complex, and calculate its cohomology. 
Taking the canonical truncation of
 this complex in the degrees $[0, 2n]$ we get a complex calculating 
the absolute Hodge cohomology of $X(\C)$. 
Let us define, following Deligne, yet another complex of abelian groups 
quasiisomorphic to the latter complex.

Let ${\cal D}_X^{ p,q} = {\cal D}^{ p,q}$ be the abelian group of complex valued distributions of type $(p,q)$ on $X(\C)$. Consider the following cohomological 
bicomplex,  where ${\cal D}_{cl}^{ n,n}$ is the subspace of closed distributions, and ${\cal D}^{0,0}$ is in degree $1$:
$$
\begin{array}{ccccccccc}
&&&&&&&&{\cal D}_{cl}^{n,n}\\
&&&&&&&&\\
&&&&&&&2 \overline \partial \partial\nearrow&\\
&&&&&&&&\\
{\cal D}^{0,n-1}&\stackrel{\partial}{\longrightarrow}&{\cal D}^{1,n-1}&\stackrel{\partial}{\longrightarrow}&
...&\stackrel{\partial}{\longrightarrow}&{\cal D}^{n-1,n-1}&&\\
&&&&&&&&\\
\overline \partial \uparrow &&\overline \partial \uparrow &&&&\overline \partial \uparrow&&\\
 ...&...&...&...&...&...&...&&\\
 \overline \partial \uparrow&&\overline \partial \uparrow&&&&\overline \partial \uparrow&&\\
&&&&&&&&\\
{\cal D}^{0,1}&\stackrel{\partial}{\longrightarrow}& {\cal D}^{1,1}&\stackrel{\partial}{\longrightarrow}&... &\stackrel{\partial}{\longrightarrow}& {\cal D}^{n-1,1}&&\\
&&&&&&&&\\
\overline \partial \uparrow&&\overline \partial \uparrow&&&&\overline \partial \uparrow&&\\
&&&&&&&&\\
{\cal D}^{0,0}&\stackrel{\partial}{\longrightarrow} &{\cal D}^{1,0}& \stackrel{\partial}{\longrightarrow} &...&\stackrel{\partial}{\longrightarrow} &{\cal D}^{n-1,0}&&
\end{array}
$$
The   complex $C^{\bullet}_{{\cal D}}(X; n)$ is a  subcomplex 
of the total complex 
of this bicomplex provided by the  $\R(n-1)$-valued
 distributions in  the $n \times n$ square of the diagram and the subspace 
   ${\cal D}_{\R, cl}^{n,n}(n) \subset {\cal D}_{cl}^{n,n}$  of
 the $\R(n)$-valued distributions of type $(n,n)$. Notice that 
 $\overline \partial \partial$ sends $\R(n-1)$-valued distributions 
to $\R(n)$-valued distributions. The cohomology of this complex of abelian
groups is the absolute Hodge cohomology of $X(\C)$, see Proposition 2.1 of 
\cite{G7}. 
Now if $X$ is a variety over $\R$, then
$$
C^{\bullet}_{{\cal D}}(X_{/\R}; n) := 
C^{\bullet}_{{\cal D}}(X; n)^{\overline F_{{\infty}}}; \quad 
H_{{\cal D}}^i(X_{/\R}, \R(n)) = H^i  C^{\bullet}_{{\cal D}}(X_{/\R};n)
$$
where $\overline F_{{\infty}}$ is the composition  of the 
involution $F_{{\infty}}$   on $X(\C)$ induced by the complex 
conjugation  with the complex conjugation of coefficients.

{\bf 3. The regulator map}. 

\begin{theorem-construction} \label{6.11.02.1} 
Let $X$ be a regular complex projective variety. Then 
there exists {\rm canonical} homomorphism of complexes 
$$
{\cal P}^{\bullet}(n): {\cal Z}^{\bullet}(X; n) 
\longrightarrow C^{\bullet}_{{\cal D}}(X; n)
$$
If $X$ is defined over $\R$ then its image   lies in  
the subcomplex $C^{\bullet}_{{\cal D}}(X_{/\R}; n)$. 
\end{theorem-construction}

To define this homomorphism we need the following construction. 
Let $X$ be a variety over $\C$ and $f_1,...,f_m$ be $m$ rational functions on $X$.
 The form 
$$
\pi_m \Bigl(
d \log f_1 \wedge ... \wedge d \log f_m\Bigr),  
$$
where $\pi_n(a + ib) =a$ if $n$ odd, and $\pi_n(a + ib) = ib$ if $n$ even, has
zero periods. It has a canonical primitive defined as follows. 
Consider the following  $(m-1)$-form on $X(\C)$:   
\begin {equation} \label{1}
r_{m-1}(f_1\wedge ...\wedge  f_m) :=
\end {equation}
$$
 {\rm Alt}_m \sum_{j\geq 0} c_{j,m}\log|f_1|d\log|f_2|
\wedge ... \wedge d\log|f_{2j+1}|\wedge di\arg f_{2j+2}\wedge ... \wedge
di\arg f_{m}
$$
Here $c_{j,m}:= ((2j+1)!(m-2j-1)!)^{-1} $ and ${\rm Alt}_m$ is the operation of alternation: 
$$
{\rm
Alt}_m F(x_1,...,x_m):= \sum_{\sigma \in S_m}(-1)^{|\sigma|}F(x_{\sigma
(1)},...,x_{\sigma (m)})
$$
So $r_{m-1}(f_1\wedge ...\wedge f_m)$ is an $\R(m-1)$-valued $(m-1)$-form. One has 
$$
d r_{m-1}(f_1\wedge ...\wedge f_m) = \pi_m \Bigl(
d \log f_1 \wedge ... \wedge d \log f_m\Bigr)  
$$ 
It is sometimes convenient to  write the form (\ref{1}) as a multiple of
$$
{\rm Alt}_m\sum_{i=1}^m(-1)^i \log |f_1| d\log f_2 \wedge ... \wedge d\log f_i 
\wedge d\log \overline f_{i+1}\wedge ... \wedge d\log \overline f_m
$$
Precisely, let 
${\cal A}^i(M)$ be the space of smooth $i$-forms on a real smooth manifold $M$. 
Consider the following map
\begin{equation} \label{6.16.04.1}
\omega_{m-1}: \Lambda^{m}{\cal A}^0(M) \lra {\cal A}^{m-1}(M)
\end{equation}
$$
\omega_{m-1}(\varphi_1 \wedge ... \wedge \varphi_{m}) := 
$$
$$
\frac{1}{m!}{\rm Alt}_{m}
\Bigl(\sum_{k=1}^{m}
(-1)^{k-1}\varphi_1\partial \varphi_2 \wedge ... \partial \varphi_k \wedge \overline 
\partial \varphi_{k+1} \wedge ... \wedge \overline \partial \varphi_{m}\Bigr)
$$
For example 
$$
\omega_0(\varphi_1) = \varphi_1;\qquad 
\omega_1(\varphi_1 \wedge \varphi_{2}) = \frac{1}{2}\Bigl(\varphi_1\partial \varphi_2  - 
\varphi_2\partial \varphi_1 - \varphi_1\overline \partial \varphi_2  +
\varphi_2\overline \partial \varphi_1   \Bigr)
$$
Then one easily checks that 
\begin{equation} \label{6.16.04.2}
d\omega_{m-1}(\varphi_1 \wedge ... \wedge \varphi_{m}) = 
\partial \varphi_1 \wedge ... \wedge \partial \varphi_{m} 
+(-1)^m \overline \partial \varphi_1 \wedge ... \wedge \overline \partial \varphi_{m} +
\end{equation}
$$
\sum_{i=1}^{m}(-1)^{i}\overline \partial \partial \varphi_i \wedge 
\omega_{m-2}(\varphi_1 \wedge ... \wedge \widehat \varphi_i \wedge ... \wedge \varphi_{m})
$$
Now let $f_i$ be rational functions on 
a complex algebraic variety $X$. Set $M:= X^0(\C)$, where $X^0$ is the open part 
of $X$ where the functions $f_i$ are regular. Then $\varphi_i := \log|f_i|$ 
are smooth functions on $M$, and we have 
$$
\omega_{m-1}(\log |f_1| \wedge ... \wedge \log |f_{m}|) = r_{m-1}(f_1 \wedge ... \wedge f_{m})
$$

Denote by ${\cal D}^*_{X, \R}(k) = {\cal D}^*_{\R}(k)$ the subspace of $\R(k)$--valued distributions
 in ${\cal D}_{X(\C)}^*$. 

Let  $Y^0$ be the nonsingular part of $Y$, and $i_Y^0: Y^0(\C) \hookrightarrow
 Y(\C)$ 
the canonical embedding. 

\begin{proposition} \label{8.6.02.2}
 Let $Y$ be an arbitrary irreducible subvariety of a smooth complex variety 
$X$ and $f_1, ..., f_m \in {\cal O}^*(Y)$. Then for any smooth differential 
form $\omega$
with compact support on $X(\C)$ the following integral is convergent:
$$
\int_{Y^0(\C)}r_{m-1}(f_1\wedge ...\wedge f_m) \wedge i_Y^0\omega
$$ 
Thus there is a distribution 
$r_{m-1}(f_1\wedge ... \wedge f_m)\delta_Y$ on 
$X(\C)$: 
$$
<r_{m-1}(f_1\wedge  ... \wedge f_m)\delta_Y, \omega>:= 
\int_{Y^0(\C)}r_{m-1}(f_1\wedge  ... \wedge f_m) \wedge i_Y^0\omega
$$
 It provides  a group homomorphism 
\begin{equation} \label{2}
r_{m-1}: \Lambda^m\C(Y)^{\ast} \longrightarrow {\cal D}^{m-1}_{X, \R}(m-1) 
\end{equation} 
\end{proposition}

{\bf Construction}. We have to construct a morphism of complexes
$$
\begin{array}{ccccccccc}
... &\longrightarrow&{\cal Z}^{1}(X; n)&\longrightarrow  &...&\longrightarrow&{\cal Z}^{2n-1}(X; n)&\longrightarrow&{\cal Z}^{2n}(X; n)\\
&&&&&&&&\\
&&\downarrow {\cal P}^{1}(n)&&... &&\downarrow {\cal P}^{2n-1}(n)&&\downarrow 
{\cal P}^{2n}(n)\\
&&&&&&&&\\
0&\lra &{\cal D}^{0,0}_{\R}(n-1)&  \stackrel{}{\longrightarrow}&...&\stackrel{}{\longrightarrow}&{\cal D}_{\R}^{n-1,n-1}(n-1)&\stackrel{2 \overline \partial \partial}{\longrightarrow}&{\cal D}_{\R}^{n,n}(n)
\end{array}
$$
Here at the bottom stays the complex ${\cal C}^{\bullet}_{\cal D}(X; n)$. 

Let $Y \in {\cal Z}^{2n}(X; n)$ 
be a codimension $n$ cycle in $X$. By definition 
$$
{\cal P}^{2n}(n)(Y):= (2\pi i)^n\delta_Y
$$

Let us construct homomorphisms
$$
{\cal P}^{2n-k}(n):  {\cal Z}^{2n-k}(X; n) 
\longrightarrow {\cal D}^{2n-k-1}_{X}, \quad k>0 
$$

Denote by  $\pi_{{\AAA}^k}$ (resp. $\pi_{X}$) the projection of 
$X \times \AAA^k$ to $\AAA^k$ (resp. $X$), and by 
 $\overline \pi_{{\AAA}^k}$ (resp. $\overline \pi_{X}$) the projection of 
$X(\C) \times \C{\Bbb P}^k$ to $\C{\Bbb P}^k$ (resp. $X(\C)$). 

The pair $(L, H)$ in $\PP^k$ defines uniquely homogeneous coordinates
$(z_0:...:z_k)$ in $\PP^k$ such that the hyperplane 
$L_i $ is given by equation $\{z_i = 0\}$ and 
the hyperplane $H$ is  $\{\sum_{i=1}^{k}
z_i =z_0\}$. Then there is an element 
\begin{equation} \label{mint3}
\frac{z_1}{z_0} \wedge... \wedge \frac{z_k }{z_0}
\in \Lambda^{k-1}\C(\AAA^k)^*
\end{equation}
Let 
$Y \in  {\cal Z}^{2n-k}(X; n)$. Restricting 
to $Y$ the inverse image of  the element 
(\ref{mint3}) by  $\pi^*_{\AAA^k}$ we get an element 
\begin{equation} \label{mint2}
g_1 \wedge... \wedge g_k \in \Lambda^{k}\C(Y)^*
\end{equation}

{\bf Remark}. This works if and only if 
the cycle $Y$ intersects  properly all codimension one faces 
of $X \times L$. Indeed, if $Y$ does not intersect properly 
one of the faces, then the equation of this face restricts to zero 
to $Y$, and so (\ref{mint2}) does not make sense. 

The element (\ref{mint2}) provides, by 
Proposition \ref{8.6.02.2}, a distribution on $X(\C) \times \C {\Bbb P}^k$. 
Pushing it down by $(2\pi i)^{n-k}\cdot\overline \pi_{X}$ we get 
the distribution $ {\cal P}^{2n-k}(n)(Y)$:

\begin{definition} \label{mint6}
$
{\cal P}^{2n-k}(n)(Y) := (2\pi i)^{n-k}\cdot
\overline \pi_{X *}r_{k-1}(g_1 \wedge... \wedge g_k) 
$. 
\end{definition}
In other words,  
the following distribution makes sense: 
$$
{\cal P}^{2n-k}(n)(Y)= \quad (2\pi i)^{n-k}
{\overline \pi_{X}}_*\left(\delta_Y \wedge \overline \pi_{\AAA^k}^{\ast}r_{k-1}\left(\frac{z_1}{z_0} \wedge... \wedge \frac{z_k }{z_0}\right)\right)
$$

One can  rewrite definition \ref{mint6} more explicitly 
as an  integral over $Y(\C)$. Namely,  
let $\omega$ be a smooth form on $X(\C)$ and  $i_{Y}: Y 
\hookrightarrow X \times  {\Bbb P}^k$.  Then 
$$
<{\cal P}^{2n-k}(n)(Y),\omega>= 
(2\pi i)^{n-k}\int_{Y^0(\C)}  
r_{k-1}(g_1\wedge ... \wedge g_k)\wedge i_{Y^0(\C)}^*\overline \pi^*_X\omega 
$$

It is easy to check that ${\cal P}^{2n-k}(n)(Y)$ lies 
in ${\cal C}^{2n-k}_{\cal D}(X; n)$. 
Therefore we defined the maps  ${\cal P}^{k}(n)$. It was proved in Theorem 2.12 in 
\cite{G7} that 
${\cal P}^{\bullet}(n)$ is a homomorphism of complexes.

{\bf 4. The Higher Arakelov Chow groups}. 
Let $X$ be a regular complex variety. 
Denote by  $\widetilde  C^{\bullet}_{{\cal D}}(X; n)$ 
the quotient of the complex 
$C^{\bullet}_{{\cal D}}(X; n)$ along the subgroup 
${\cal A}_{cl, \R}^{n,n}(n) \subset {\cal D}_{cl, \R}^{n,n}(n) $ of 
closed smooth forms. 
The cone of the homomorphism ${\cal P}^{\bullet}(n)$ shifted by $-1$ 
is the Arakelov motivic complex: \index{Arakelov motivic complex}:
$$
\widehat  {\cal Z}^\bullet(X; n):= {\rm Cone}
\Bigl( {\cal Z}^\bullet(X; n)
 \stackrel{ }{\longrightarrow} \widetilde  C^{\bullet}_{{\cal D}}(X; n)
\Bigr)[-1]
$$
The Higher Arakelov Chow groups \index{Higher Arakelov Chow groups} are its cohomology:
$$
\widehat  {CH}^n(X; i):= H^{2n-i}(\widehat  {\cal Z}^{\bullet}(X; n))
$$
 Recall the arithmetic Chow groups  defined by Gillet-Soul\'e 
\cite{GS} as follows:
$$
\widehat  {CH}^n(X):= 
$$
\begin{equation} \label{6.11.02.10}
\frac{\{(Z,g); \frac { \overline \partial \partial}{\pi i} 
g+\delta_Z \in {\cal A}^{n,n}\}}{\{(0,\partial u + \overline \partial v); 
({\rm div} f, -\log|f|), f \in \C(Y), 
{\rm codim} (Y) = n-1\}}
\end{equation} 
Here $Z$ is a divisor in $X$, $f$ is a rational function on a divisor $Y$ in $X$, 
$$
g \in {\cal D}_{\R}^{n-1,n-1}(n-1), \quad (u,v) \in C^{2n-2}_{{\cal D}}(X; n) = 
({\cal D}^{n-2,n-1} \oplus {\cal D}^{n-1,n-2})_{\R}(n-1)
$$

\begin{proposition} \label{ch}
 $\widehat  {CH}^n(X; 0) = \widehat  {CH}^n(X)$.
\end{proposition}

{\bf Proof}.  Let us look at the very right part of the complex 
$\widehat  {\cal Z}^{\bullet}(X; n)$:
$$
\begin{array}{ccccc}
...&\longrightarrow&{\cal Z}^{2n-1}(X; n)&\longrightarrow&{\cal Z}^{2n}(X; n) \\
&&&&\\
&&\downarrow {\cal P}^{2n-1}(n)&& \downarrow {\cal P}^{2n}(n)\\
&&&&\\
({\cal D}^{n-2,n-1} \oplus {\cal D}^{n-1,n-2})_{\R}(n-1)&
\stackrel{(\partial , \overline \partial)}{\longrightarrow}&
{\cal D}_{\R}^{n-1,n-1}(n-1)&\stackrel{2\overline \partial \partial}
{\longrightarrow}&{\cal D}_{\R}^{n,n}(n)/{\cal A}_{\R}^{n,n}(n)
\end{array}
$$
Consider the stupid truncation of the  Gersten complex on $X$:
\begin{equation} \label{3.12.03.1}
\prod_{Y \in X_{n-1}}\C(Y)^* 
\longrightarrow {\cal Z}_0(X; n) 
\end{equation}
 It maps to the stupid truncation  
$\sigma_{\geq 2n-1}\widehat  {\cal Z}^{\bullet}(X; n)$  of the cycle  complex
as follows. The isomorphism ${\cal Z}_0(X; n) = {\cal Z}^{2n}(X; n)$ provides 
the right component of the map. 
A pair $(Y;f)$ where $Y$ is an irreducible codimension $n-1$ subvariety of $X$ 
maps to the cycle $(y,f(y)) \subset X \times (\PP^1 - \{1\})$. 
It is well known that such cycles $(y,f(y))$ plus 
$\partial \widehat  {\cal Z}^{2n-2}(X; n)$ 
generate $\widehat  {\cal Z}^{2n-1}(X; n)$. 
 Computing the composition of this map with the homomorphism 
${\cal P}^{\bullet}(n)$ we end up precisely with 
formula (\ref{6.11.02.10}). 
The proposition is proved.

{\bf 5. Remarks on special values of Dedekind $\zeta$-functions}. 
Let $\widetilde {\Gamma}_{\cal A}({\cal O}_F, 1)$ be
the three term complex (\ref{3.21.03.5}). It 
consists of locally compact abelian groups.
Each of them is equipped with a natural Haar measure. Indeed, 
the measure of  a discrete group is normalized so that 
the measure of the identity element is $1$;   
the group $\R$ has the canonical measure $dx$; and we use the 
product measure for the products. 
We need the following general observation. 

\begin{lemma-definition} \label{4.10.03.1} Let 
$$
A^{\bullet} = \quad ... \lra  
A_1 \lra A_2 \lra A_3 \lra ... 
$$
 be a  
complex of locally compact abelian groups 
such that  

i) Each of the groups $A_i$ is equipped with an invariant 
Haar measure $\mu_i$. 

ii) The cohomology groups are compact. 

iii) Only finite number of the cohomology groups are nontrivial, 
and almost all groups $A_i$ are discrete groups with canonical measures.

Then  
there is a naturally defined number
${\rm R}_{\mu}A^{\bullet}$, and  ${\rm R}_{\mu}A^{\bullet}[1] = 
({\rm R}_{\mu}A^{\bullet})^{-1}$.  
\end{lemma-definition}

{\bf Construction}. Let 
\begin{equation} \label{3.22.03.1}
0 \lra A \lra B \lra C \lra 0
\end{equation} 
 be an exact sequence of 
locally compact abelian groups. Then a choice of Haar measure 
for any two of the groups $A, B, C$ determines  naturally 
the third one. For example Haar measures 
$\mu_A, \mu_C$ on $A$ and $C$ determine the following  Haar measure
$\mu_{A,C}$ on $B$. Take a compact 
subset of $C$ and its section $K_C \subset B$, and take a compact 
subset $K_A$. Then $\mu_{A, C}(K_A \cdot K_C) := \mu_A(K_A) \mu_C(K_C)$.

For the complex (\ref{3.22.03.1}) placed in degrees $[0,2]$ 
we put   ${\rm R}_{\mu}(\ref{3.22.03.1}) = \mu_{A,C}/\mu_B$.

Let us treat first  the case when $A^{\bullet}$ is a finite complex. 
Then we define the invariant ${\rm R}_{\mu}A^{\bullet}$ by induction. 
If $A^{\bullet}= A[-i]$ is concentrated in just one degree, $i$, that 
$A$ is compact and equipped with Haar measure $\mu$. We put 
${\rm R}_{\mu}A^{\bullet}:= \mu(A)^{(-1)^i}$. 
Assume $A^{\bullet}$ starts from $A_0$. One has 
$$
{\rm Ker}f_0 \lra A_0 \lra {\rm Im}f_0 \lra 0; \qquad  
{\rm Im}f_0 \lra A_1 \lra  A_1/{\rm Im}f_0\lra 0 
$$
Since ${\rm Ker}f_0$ is compact, we can choose the 
volume one Haar measure on it. 
This measure and the measure 
$\mu_0$ on $A_0$ provides, via the first short exact sequence, 
a measure on ${\rm Im}f_0$. Similarly using this and the second exact sequence 
we get a measure on $A_1/{\rm Im}f_0$. 
Therefore we have the measures on the truncated complex 
$\tau_{\geq 1}A^{\bullet}$. 
Now we define 
$$
{\rm R}_{\mu}A^{\bullet}:= \mu_0(A_0)\cdot {\rm R}_{\mu}\tau_{\geq 1}A^{\bullet}
$$
If $A^{\bullet}$ is an infinite complex 
we set ${\rm R}_{\mu}A^{\bullet}:= {\rm R}_{\mu}(\tau_{-N,N}A^{\bullet})$ 
for sufficiently big $N$.  
Here $\tau_{-N,N}$ is the canonical truncation functor. Thanks to iii) this 
 does not depend on the choice of $N$. 
The lemma is proved. 

Now the class number formula (\ref{3.21.03.1}) reads
\begin{equation} \label{3.03.21.6q}
\lim_{s \to 0}s^{-(r_1+r_2-1)}\zeta_F(s) = -
{\rm R}_{\mu}\widetilde {\Gamma}_{\cal A}({\cal O}_F, 1)
\end{equation}

Lichtenbaum's conjectures \cite{Li1}
on the special values of the Dedekind $\zeta$-functions 
can be reformulated in a similar way: 
$$
\zeta_F(1-n) \stackrel{?}{=} \pm {\rm R}_{\mu}{\Gamma}_{\cal A}({\cal O}_F, n); 
\quad n>1
 $$
It is not quite clear what is the most natural normalization of the 
regulator map. 
In the classical $n=1$ case this 
formula needs modification, as was explained in s. 1.1, to take into account 
the pole of the $\zeta$--function. 

{\bf Example}. The group 
$H^2\Gamma_{\cal A}({\cal O}_F, 2)$ sits in the exact sequence
$$
0 \lra R_2(F) \lra H^2\Gamma({\cal O}_F, 2) \lra K_2({\cal O}_F) \lra 0 
$$
To calculate it let us define the Bloch-Suslin complex for ${\rm Spec}({\cal O}_F)$:
\begin{equation} \label{4.12.03.2}
B({\cal O}_F,2): \qquad B_2(F) \lra \Lambda^2F^* \lra \prod_{\cal P}k^*_{\cal P}
\end{equation}
 Its Arakelov version is the total complex of the following bicomplex,
where the vertical map is given by the dilogarithm: $\{x\}_2 
\lms ({\cal L}_2(\sigma_1(x), ..., {\cal L}_2(\sigma_{r_2}(x))$. 
\begin{equation} \label{4.12.03.1}
\begin{array}{ccccc}
\R^{r_2}&&&&\\
\uparrow &&&&\\
B_2(F)& \lra &\Lambda^2F^* & \lra&\prod_{\cal P}k^*_{\cal P}
\end{array}
\end{equation}
Then $H^2\Gamma_{\cal A}({\cal O}_F, 2) = H^2B_{\cal A}({\cal O}_F, 2)$.
The second map in (\ref{4.12.03.2}) 
 is surjective (\cite{Mi}, cor. 16.2). 
Using this one can check that 
$$
H^i\Gamma_{\cal A}({\cal O}_F, 2) = 0\quad \mbox{for $i\geq 3$}
$$ 
(Notice that  $H^1B_{\cal A}({\cal O}_F, 2) \not = 
 H^1\Gamma_{\cal A}({\cal O}_F, 2)$). Summarizing, we should have 
$$
\zeta_F(-1) \stackrel{?}{=} \pm {\rm R}_{\mu}{\Gamma}_{\cal A}({\cal O}_F, 2) 
= \pm \frac{{\rm vol} H^2{\Gamma}_{\cal A}({\cal O}_F, 2)}{
|H^1{\Gamma}_{\cal A}({\cal O}_F, 2)|}
$$
For totally real fields it is a version of the Birch-Tate conjecture. 
It would be very interesting to compare this with the 
Bloch--Kato conjecture.

\section{Grassmannian polylogarithms and Borel's regulator}

{\bf 1. The Grassmannian polylogarithm \cite{G4}}. \index{Grassmannian polylogarithms} 
Let $h_1,...,h_{2n}$ be arbitrary $2n$ hyperplanes in $\C \PP^{n-1}$. Choose an additional hyperplane $h_0$. Let $f_i$ be a
rational function on $\C \PP^{n-1}$ with divisor $h_i - h_0$. It is 
defined up to a scalar
factor. 
Set
$$
{\cal L}^G_n(h_1,...,h_{2n}):= (2\pi i)^{1-n}\int_{\C \PP^{n-1}}r_{2n-2}
(\sum_{j=1}^{2n}(-1)^j f_1\wedge ...
\wedge \widehat  f_j \wedge ... \wedge f_{2n})
$$
It is skew-symmetric by definition. It is easy to see that it 
does not depend on the choice of scalar in the definition of $f_i$. 
To check that it does not depend on the choice of 
$h_0$ observe that
$$
\sum_{j=1}^{2n}(-1)^j f_1\wedge ...
\wedge \widehat  f_j \wedge ... \wedge f_{2n} = \frac{f_1}{f_{2n}} \wedge
\frac{f_2}{f_{2n}} \wedge ... \wedge \frac{f_{2n-1}}{f_{2n}} 
$$
So if we  choose rational
functions $g_1,...,g_{2n-1}$ such that ${\rm div} g_i = h_i - h_{2n}$ then
$$
{\cal L}^G_n(h_1,...,h_{2n}) = (2\pi i)^{1-n}
\int_{\C \PP^{n-1}}r_{2n-2}(g_1\wedge ... \wedge g_{2n-1})
$$
{\bf Remark}. The function ${\cal L}^G_n$ is defined on the set 
of all configurations of $2n$ hyperplanes in $\C\PP^{n-1}$. However it is 
not even continuous on this set. It is real 
analytic on the submanifold of generic 
configurations.  

\begin {theorem}  
\label {0.3} The function ${\cal L}^G_n$ satisfies the following 
functional equations:

a) For any $2n+1$ hyperplanes $h_1, ..., h_{2n+1}$ in $\C \PP^{n}$ one has
\begin{equation} \label{4}
 \sum_{j=1}^{2n+1}(-1)^j {\cal L}^G_n(h_j \cap h_1,..., h_j \cap h_{2n+1}) =0
\end{equation}

b) For any $2n+1$ hyperplanes $h_1, ..., h_{2n+1}$ in $\C \PP^{n-1}$ one has
\begin{equation} \label{5}
 \sum_{j=1}^{2n+1}(-1)^j {\cal L}^G_n(h_1,...,\widehat  h_j,...,h_{2n+1})=0
\end{equation}
\end {theorem}

{\bf Proof}. a) Let $f_1,...,f_{2n+1}$ be rational functions on $\C \PP^{n}$ 
as above. 
Then
\begin{equation} \label{8}
 dr_{2n-1}\Bigl(\sum_{j=1}^{2n+1}(-1)^j f_1 \wedge ...
\wedge \widehat  f_j \wedge ... f_{2n+1}\Bigr) = 
\end{equation}
$$
\sum_{j \not = i}(-1)^{j+i-1} 2 \pi i \delta(f_j)  r_{2n-2}\Bigl(
f_1\wedge ... \widehat  f_i \wedge ...
\wedge \widehat  f_j \wedge ... f_{2n+1}\Bigr)
$$
(Notice that $d\log f_1 \wedge ... \wedge \widehat  {d\log f_j} \wedge ...
\wedge d\log f_{2n+1} =0$ on $\C \PP^{n}$). 
Integrating (\ref{8}) over $\C \PP^{n}$ we get a).

b) is obvious: we apply $r_{2n-1}$ to zero element.
Theorem is proved.

There is a  natural bijection
$$
PG^{n-1}_n:= \mbox{ $\{ (n-1)$--planes in 
 $\PP^{2n-1}$ in generic position to a simplex $L$\}}/({\Bbb G}_m^*)^{2n-1} \quad 
$$
$$  <--> \quad \left \{ \mbox{Configurations of} \quad 2n \quad\mbox {generic hyperplanes in } \PP^{n-1}\right \}
$$
given by intersecting of an $(n-1)$--plane $h$ with the codimension one faces of $L$. 

\begin{center}
\hspace{4.0cm}
\epsffile{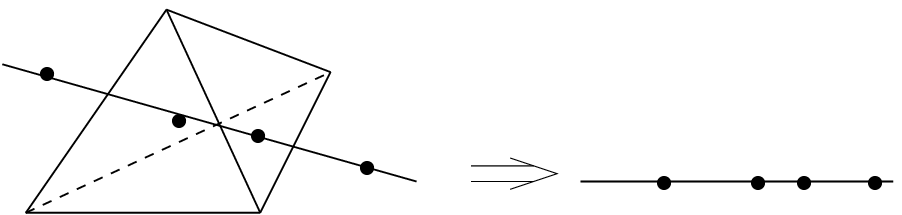}
\end{center}

Thus ${\cal L}^G_n$ is a function on the torus quotient  $PG^{n-1}_n$ of 
the generic part of Grassmannian.  


Gelfand and MacPherson \cite{GM} suggested 
a beautiful construction of the real valued version of 
$2n$-logarithms on $PG^{2n-1}_{2n}(\R)$. The construction uses the Pontryagin form. 
These functions generalize 
the Rogers dilogarithm.   

The defined above functions ${\cal L}^G_n$ on complex Grassmannians 
generalize the Bloch-Wigner dilogarithm.  They are related to the Chern classes. 
It would be very interesting to find a link between the construction 
in \cite{GM} with our construction. 

The existence of the multivalued analytic Grassmannian $n$-logarithms 
on complex Grassmannians was conjectured 
in \cite{BMS}. 
 They were constructed  in \cite{HM1}-\cite{HM2} and,  
as a particular case of the analytic Chow polylogarithms, in \cite{G4}.

Recall the following general construction. 
Let $X$ be a $G$-set and $F$ a $G$-invariant function on $X^n$ satisfying 
$$
\sum_{i=1}^n(-1)^{i}F(x_1, ..., \widehat x_i, ..., x_n) =0
$$
 Choose a point $x \in X$. Then there is an $(n-1)$-cocycle of the group $G$:
$$
f_x(g_1, ..., g_n):= F(g_1 x, ..., g_n x)
$$  

\begin{lemma} \label{point}
The cohomology class of the cocycle $f_x$ does not depend on $x$.
\end{lemma}

Thus thanks to (\ref{5}) the function ${\cal L}^G_n$ 
provides a measurable cocycle of $GL_n(\C)$. 
We want to determine its cohomology class, but a priori it is not 
even clear that it is non zero. 
To handle this problem we will show below that 
the function ${\cal L}^G_n$ is a boundary value of a certain function 
$\psi_n$ 
defined on the configurations of $2n$ points inside of the 
symmetric space $SL_m(\C)/SU(n)$. The cohomology class 
of $SL_n(\C)$ provided by this function is obviously 
related to the so-called Borel class.  
Using this we will show that the Grassmannian $n$-logarithm 
function ${\cal L}^G_n$ provides the Borel class, and moreover 
can be used to define the Borel regulator. 

Finally, the restriction of the function ${\cal L}^G_n$ 
to certain special stratum in the configuration space of $2n$ hyperplanes 
in $\C\PP^{n-1}$ provides a single valued version 
of the classical $n$--logarithm function, see sections 4-5 below. 

{\bf 2.  The function $\psi_n$.} Let $V_n$ be an $n$-dimensional complex vector space. Let
$$
\HH_n: = \left \{ \mbox{ positive definite Hermitian forms in  
   $V_n$ }\right \}/\R ^*_+ = SL_n(\C)/SU(n) 
$$
$$
= \left \{ \mbox{ positive definite Hermitian forms in
   $V_n$ with determinant } =1 \right \} 
$$
It is a symmetric space    of rank $n-1$.  For example 
$\HH_2 = {\cal H}_3$ is the hyperbolic 3-space. 
Replacing positive definite by non negative definite Hermitian forms 
we get a compactification $\overline \HH_n$ of the symmetric space $\HH_n$.

Let $G_x$ be the subgroup of $SL_N(\C)$ stabilizing the point  $x \in \HH_n$ . 
A point $x$ defines a one dimensional vector space $M_x$:
$$
x \in \HH_n \longmapsto M_x:= \left \{ \mbox{measures on    } \C \PP^{n-1} 
\mbox{ invariant under }  G_x\right \}
$$
   Namely, a point $x$ corresponds to a hermitian metric in $V_n$.  
This metric provides the Fubini-Studi Kahler form on $\C \PP^{n-1} = P(V_n)$. 
Its imaginary part is a symplectic form. Raising it to $(n-1)$-th power we get 
the Fubini-Studi volume form.
The elements of   $M_x$ are its multiples.

Let $x_0,...,x_{2n-1}$ be points  of the symmetric space $SL_n(\C)/SU(n)$.
Consider the following function
$$
\psi_n(x_0,...,x_{2n-1}) := \int_{\C \PP^{n-1}}
 \log | \frac{\mu_{x_1}}{\mu_{x_0}}| d\log|\frac{\mu_{x_2}}{\mu_{x_0}} | \wedge ... \wedge d\log| \frac{\mu_{x_{2n-1}}}{\mu_{x_0}} |
$$

  More generally, let $X$ be an $m$-dimensional manifold. 
For any $m+2$ measures $\mu_0,...,\mu_{m+1}$
 on $ X$ such that $\frac{\mu_{i}}{\mu_j}$ are smooth functions 
consider the following differential $m$-form  on $X$:
$$
\overline r_{ m}(\mu_0:...:\mu_{ m+1}) := \log | \frac{\mu_1}{\mu_0}| d\log|\frac{\mu_2}{\mu_0} | \wedge ... \wedge d\log| \frac{\mu_{ m+1}}{\mu_0} |
$$

\begin{proposition} \label{1.2.}
The integral 
$
\int_{ X}\overline r_{ m}(\mu_0:...:\mu_{ m+1}) 
$
 satisfies the following properties:

1) Skew symmetry with respect to the permutations of $\mu_i$.

2) Homogeneity: if $\lambda_i \in \R^*$ then 
$$
\int_{ X}\overline r_{ m}(\lambda_0  \mu_0: ... :\lambda_{ m+1}   \mu_{ m+1}) = \int_{X }\overline r_{ m}(\mu_0:...:\mu_{ m+1}) 
$$

3)Additivity: for any $ m+3$ measures $\mu_i$ on $X $ one has
$$
  \sum_{i=0}^{m+2} (-1)^i \int_{X }\overline r_{ m}(\mu_0: ... :\widehat \mu_i: ... :\mu_{ m+2}) =0
$$

4)  Let $g$ be a diffeomorphism of $X$.  Then   
$$
\int_{X }\overline r_{ m}(g^*\mu_0: ... :g^*\mu_{ m+1})  = \int_{X }\overline r_{ m}( \mu_0: ... : \mu_{ m+1})
$$
\end{proposition}

{\bf 3. The Grassmannian polylogarithm ${\cal L}_n^G$ is the boundary value
 of the function $\psi_n$ \cite{G7}}.  
Let $(z_0:...:z_{n-1})$ be   homogeneous coordinates in $\PP^{n-1}$. 
Let
$$
\sigma_n(z, dz):= \sum_{i=0}^{n-1} (-1)^i z_i dz_0 \wedge ... \wedge {\widehat 
{dz_i}} \wedge ... \wedge dz_{n-1}  
$$
\begin{equation} \label{FUBS}
\omega_{FS}(H) := \frac{1}{(2\pi i)^{n-1}}\frac{\sigma_n(z,dz) 
\wedge \sigma_n(\overline z, d \overline {z})}{H(z, \overline z)^{n}}
\end{equation} 
This form is clearly invariant under the group preserving the 
Hermitian form $H$. In fact  
it is the Fubini-Studi volume form.

Take any $2n$ non zero nonnegative definite Hermitian forms 
$H_0, ..., H_{2n-1}$, possibly degenerate. For each of the forms $H_i$ 
 choose a multiple  $\mu_{H_i}$ of the Fubini-Studi form 
given by formula (\ref{FUBS}). It is a volume form 
with singularities along the projectivization of 
kernel of $H_i$. 

\begin{lemma} \label{222222}
The following integral is convergent
$$
\psi_n(H_0,...,H_{2n-1}) := \int_{\C \PP^{n-1}}
 \log | \frac{\mu_{H_1}}{\mu_{H_0}}| d\log|\frac{\mu_{H_2}}{\mu_{H_0}} | \wedge ... \wedge 
d\log| \frac{\mu_{H_{2n-1}}}{\mu_{H_0}} |  
$$
\end{lemma}

This integral does not change if we multiply one of the 
Hermitian forms by a positive scalar.  
Therefore we can  extend  $\psi_n$ to a  function 
on the configuration space of $2n$ points in the compactification 
$\overline \HH_{n-1}$. This function is discontinuous.

One can realize 
$\C \PP^{n-1}$ as the smallest stratum  of the boundary of $\HH_n$. Indeed, 
let $\C\PP^{n-1} = P(V_n)$. For 
a hyperplane $h \in V_n$ let 
$$
F_h:= \left \{ \mbox{    nonnegative definite hermitian forms in 
  $V_n$ with  kernel $h$}   \right \}/ \R _+^*
$$
The set of 
 hermitian forms in $V_n$ with the 
kernel  $h$  is isomorphic to $\R _+^*$,  so  $F_h$ defines a point on 
the boundary of $\overline \HH_n$. 
Therefore  Lemma 
\ref{222222} provides a function 
\begin{equation} \label{1221N} 
\psi_n(h_0,...,h_{2n-1}) := \psi_n(F_{h_0}, ..., F_{h_{2n-1}})
\end{equation}

Applying Lemma \ref{point} to the case when $X$ is 
$\overline {\Bbb H}_n$ and using only the fact that 
the function $\psi_n(x_0,...,x_{2n-1})$ is well defined for {\it any}
 $2n$ points in $\overline {\Bbb H}_n$ and satisfies 
the cocycle condition for any $2n+1$ of them we get 

\begin{proposition} \label{pointx1}
Let $x \in {\Bbb H}_n$ and $h$ is a hyperplane in $\C\PP^{n-1}$. 
Then the cohomology classes of the following cocycles coincide:
$$
\psi_n(g_0 x,...,g_{2n-1}x) \quad \mbox{and} \quad \psi_n(g_0 h,...,g_{2n-1}h)
$$ 
\end{proposition}

The Fubini-Studi volume form corresponding to a hermitian 
form from the set $F_{h}$ is a
 Lebesgue measure on the affine space $\C\PP^{n-1} -  h$. 
Indeed, if 
$h_0 = \{z_0 =0\}$ then (\ref{FUBS}) specializes to 
    $$
\frac{1}{(2\pi i)^{n-1}}  d\frac{z_1}{z_0} \wedge ... \wedge 
d\frac{z_{n-1}}{z_0}  \wedge  d\frac{ \overline z_1}{\overline z_0} 
\wedge ... \wedge d\frac{ \overline z_{n-1}}{\overline  z_0}  
$$
Using this it is easy to prove the following proposition.

\begin{proposition} \label{1/18.1}
For any $2n$ hyperplanes $h_0,...,h_{2n-1}$ in $\C \PP^{n-1}$ 
one has 
$$
\psi_n(h_0,...,h_{2n-1}) = (-4)^{-n} \cdot (2\pi i)^{n-1} (2n)^{2n-1}{2n-2\choose n-1} 
\cdot {\cal L}_n^G(h_0, ..., h_{2n-1})
$$
\end{proposition}

{\bf 4. Construction of the Borel regulator 
via Grassmannian polylogarithms}. We start from a normalization of the Borel class $b_n$. 
Denote by $H^*_c(G, \R)$ the continuous cohomology of a Lie group $G$. 
Let us define an isomorphism
$$
\gamma_{\rm DR}: H_{\rm DR}^{k} (SL_n(\C), \Q) 
\stackrel{\sim}{\longrightarrow}
H^{k}_c(SL_n(\C), \C)
$$ 
We do it in two steps. First, let us define an isomorphism
$$
\alpha: H_{\rm DR}^{k} (SL_n(\C), \C) \stackrel{\sim}{\longrightarrow} 
{\cal A}^{k}(SL_n(\C)/SU(n))^{SL_n(\C)} \otimes \C
$$
It is well known that any cohomology class on the left is represented 
by a biinvariant, and hence closed, 
 differential $k$--form $\Omega$ on $SL_n(\C)$. Let us restrict it first 
to the Lie algebra, and then  
to the orthogonal 
complement $su(n)^{\perp}$ to the Lie subalgebra $su(n) \subset  sl_n(\C)$. 
Let $e$ be the point of  ${\Bbb H}_n$ corresponding to the  subgroup $SU(n)$. 
We identify the $\R$--vector spaces $T_e{\Bbb H}_n$ and  $su(n)^{\perp}$. 
The obtained exterior form on $T_e{\Bbb H}_n$ is restriction of 
 an invariant closed differential form $\omega$  
on the symmetric space 
${\Bbb H}_n$. 

Now let us construct, following J. Dupont \cite{Du}, an isomorphism 
$$
\beta: {\cal A}^{k}(SL_n(\C)/SU(n))^{SL_n(\C)} \stackrel{\sim}{\lra} H^{k}_c(SL_n(\C), \R)
$$
For any ordered $m+1$ points $x_1,...,x_{m+1}$ 
in ${\Bbb H}_n$ there is 
a  geodesic 
simplex $I(x_1,...,x_{m+1})$ in ${\Bbb H}_n$. It is constructed inductively 
as follows. Let 
$I(x_1,x_2)$ be the geodesic from $x_1$ to $x_2$. The  
geodesics from $x_3$ to 
the points of  $I(x_1,x_2)$ form a geodesic triangle 
$I(x_1,x_2,x_3)$, and so on. If $n>2$ 
the geodesic simplex $I(x_1,...,x_k)$ depends on the 
order of vertices. 

Let $\omega$ be an invariant 
 differential $m$-form on 
$SL_n(\C)/SU(n)$. Then it is closed, and 
 provides a volume of the geodesic simplex: 
$$
{\rm vol_{\omega}} I(x_1,...,x_{m+1}):= \int_{I(x_1,...,x_{m+1})}\omega
$$
The boundary of the simplex 
$I(x_1,...,x_{m+2})$ is the alternated sum of simplices 
$I(x_1,..., \widehat x_i,..., x_{m+2})$. Since the form $\omega$ is 
closed, the Stokes theorem yields 
$$
\sum_{i=1}^{m+2}(-1)^i \int_{I(x_1,..., \widehat x_i,..., x_{m+2})}\omega = 
\int_{I(x_1,..., x_{m+2})}d\omega = 0
$$
This just means that for a given point $x$ the function
$ {\rm vol_{\omega}} I(g_1 x,...,g_{m+1} x)$ is a smooth $m$-cocycle 
of the Lie group 
$SL_n(\C)$. 
By Lemma \ref{point} cocycles corresponding to different points $x$ 
are canonically cohomologous. 
The obtained cohomology class is the class  $\beta(\omega)$. 
Set $\gamma_{\rm DR}:= \beta \circ \alpha$. 

It is known that 
$$
H^*_{\rm DR}(SL_n(\C), \Q) = \Lambda_{\Q}^*(C_3, ..., C_{2n-1})
$$
where 
$$
C_{2n-1}:= {\rm tr}(g^{-1}dg)^{2n-1} \in \Omega^{2n-1}(SL)
$$
The Hodge considerations 
shows that $[C_{2n-1}] \in  H_{\rm Betti}^{2n-1}(SL_n(\C), \Q(n))$. 

\begin{lemma} \label{7.1.02.3} $\alpha(C_{2n-1})$ is an 
$\R(n-1)$--valued differential form. So it provides 
a cohomology  class 
$$
 b_n := \gamma_{\rm DR}(C_{2n-1}) \in H^{2n-1}_c(SL_n(\C), \R(n-1)) 
$$
\end{lemma}
We call the cohomology class provided by this lemma the Borel class, and 
use it below to construct the Borel regulator. 

It is not hard to show that the cohomology class of the cocycle 
$\psi_n(g_0x, ..., g_{2n-1}x)$ is a non zero multiple of the Borel class. 
So thanks to propositions \ref{pointx1} and \ref{1/18.1} the same is true 
for the cohomology class provided by the Grassmannian $n$--logarithm. 
The final result will be stated in Theorem \ref{6.24.02.7} below. 

Now let us construct the Borel regulator 
via Grassmannian polylogarithms. 
Let $G$ be a group. The diagonal map $\Delta: G \lra G \times G$ provides a homomorphism $\Delta_*: H_n(G) \lra H_n(G \times G)$. 
Recall that 
$$
{\rm Prim}H_nG := \{x \in H_n(G)| \Delta_*(x) = x \otimes 1 + 1 \otimes x\}
$$
Set $A_{\Q}:= A \otimes \Q$. One has 
$$
K_n(F)_{\Q} = {\rm Prim}H_nGL(F)_{\Q} = {\rm Prim}H_nGL_n(F)_{\Q}
$$
where the second isomorphism is provided by Suslin's stabilization theorem. 
Let 
$$
B_n \in H_c^{2n-1}(GL_{2n-1}(\C), \R(n-1))
$$ be a cohomology 
which goes to $b_n$ under the restriction map to $GL_n$. We define 
the Borel regulator map by restricting  the class $B_n$ to
the subspace $K_{2n-1}(\C)_{\Q}$ of $ H_{2n-1}(GL_{2n-1}(\C), {\Q})$: 
$$
r_{n}^{\rm Bo}(b_n): = <B_n, \ast>: K_{2n-1}(\C)_{\Q} \lra \R(n-1)
$$
It does not depend on the choice of $B_n$. 

Recall the Grassmannian complex $C_*(n)$
$$
... \stackrel{d}{\lra} C_{2n-1}(n) \stackrel{d}{\lra}  C_{2n-2}(n) 
\stackrel{d}{\lra} ... \stackrel{d}{\lra} C_{0}(n)
$$
where $C_k(n)$ is the free abelian group generated by configurations, i.e. 
$GL(V)$-coinvariants,  
of $k+1$ 
vectors $(l_0, ..., l_{k})$ in generic position in an $n$--dimensional 
vector space $V$ over a field $F$, and $d$ is given by the standard 
formula 
\begin{equation} \label{3.13.03.1}
(l_0, ..., l_k) \lms \sum_{i=0}^k(-1)^i(l_0, ..., \widehat l_i, ..., l_k)
\end{equation} 
The group $C_k(n)$ is in degree $k$. 
Since it is a homological resolution of the trivial $GL_n(F)$--module $\Z$  
(see Lemma 3.1 in \cite{G2}), 
there is canonical homomorphism 
$$
\varphi_{2n-1}^n: H_{2n-1}(GL_n(F)) \lra H_{2n-1}(C_*(n)) 
$$
Thanks to Lemma \ref{7.1.02.3} 
the Grassmannian $n$--logarithm function provides a homomorphism 
\begin{equation} \label{6.24.02.3}
{\cal L}_n^G: C_{2n-1}(n) \lra \R(n-1); \quad (l_0, ..., l_{2n-1}) \lms 
{\cal L}^G_n(l_0, ..., l_{2n-1})
\end{equation}
Thanks to the functional equation (\ref{4}) for ${\cal L}_n^G$ 
it is zero on the subgroup $dC_{2n}(n)$. So it induces a homomorphism
$$
{\cal L}_n^G: H_{2n-1}(C_{*}(n)) \lra \R(n-1);
$$

Let us  extend the map $\varphi_{2n-1}^n\circ {\cal L}^G_n$ 
to a homomorphism from $H_{2n-1}(GL_{2n-1}(\C))$ to $\R(n-1)$ 
following  s. 3.10 of \cite{G2},   
Consider the following bicomplex:
$$
\begin{array}{cccccccc}
&&&&&... & \stackrel{d}{\lra}&C_{2n-1}(2n-1)\\
&&& &&&&\downarrow \\
&&&...&&...&&...\\
&&& \downarrow &&&&\downarrow \\
& ... & \stackrel{d}{\lra} & C_{2n-1}(n+1) & \stackrel{d}{\lra} & ... & 
\stackrel{d}{\lra} & C_{n+1}(n)\\
&\downarrow && \downarrow &&&&\downarrow \\
... \stackrel{d}{\lra} & C_{2n-1}(n) & \stackrel{d}{\lra} & C_{2n-2}(n) & 
\stackrel{d}{\lra}& ... &\stackrel{d}{\lra} & C_n(n) 
\end{array}
$$
The horizontal differentials are given by formula 
(\ref{3.13.03.1}), and the vertical by 
$$
(l_0, ..., l_k) \lms \sum_{i=0}^k(-1)^i(l_i|l_1, ..., \widehat l_i, ..., l_k)
$$
Here $(l_i|l_1, ..., \widehat l_i, ..., l_k)$ means projection of the configuration 
$(l_1, ..., \widehat l_i, ..., l_k)$ to the quotient $V/<l_i>$. 
The total complex of this bicomplex is 
 called the weight $n$ bi--Grassmannian complex $BC_*(n)$.  

Let us  extend homomorphism 
(\ref{6.24.02.3}) to a homomorphism
$$
{\cal L}_n^G: BC_{2n-1}(n) \lra \R(n-1)
$$
by setting it zero on the groups $C_{2n-1}(n+i)$ for $i>0$. 
The functional equation  (\ref{5}) for 
the Grassmannian $n$--logarithm    
just means that the composition 
$$
C_{2n}(n+1) \lra C_{2n-1}(n) \stackrel{{\cal L}_n^G}{\lra} \R(n-1),
$$
where the first map is a vertical arrow in $BC_*(n)$, is zero. 
Therefore we get a homomorphism
$$
{\cal L}_n^G: H_{2n-1}(BC_{*}(n)) \lra \R(n-1)
$$ 
 The bottom row of the Grassmannian bicomplex is the stupid truncation of the 
Grassmannian complex at
the group $C_n(n)$. 
So there is a homomorphism 
\begin{equation} \label{6.24.02.2}
H_{2n-1}(C_*(n)) \lra H_{2n-1}(BC_*(n)) 
\end{equation} 
In \cite{G1}-\cite{G2} we proved that there are homomorphisms 
$$
\varphi_{2n-1}^m:  H_{2n-1}(GL_m(F)) \lra H_{2n-1}(BC_*(n)), 
\quad m \geq n
$$ 
whose restriction to the subgroup $GL_n(F)$ coincides 
with the composition 
$$
H_{2n-1}(GL_n(F)) \stackrel{{\varphi_{2n-1}^n}}{\lra} H_{2n-1}(C_*(n)) 
\stackrel{(\ref{6.24.02.2})}{\lra} H_{2n-1}(BC_*(n)), 
$$

\begin{theorem}  \label{6.24.02.7}
The composition 
$$
 K_{2n-1}(\C) \stackrel{\sim}{\lra} 
{\rm Prim} H_{2n-1}(GL_{2n-1}(\C), \Q) \stackrel{\varphi_{2n-1}^{2n-1}}{\lra} 
H_{2n-1}(BC_{*}(n)_{\Q}) \stackrel{{\cal L}_n^G}{\lra} \R(n-1)
$$ 
equals to 
$$
-(-1)^{n(n+1)/2} \cdot \frac{(n-1)!^2}{n(2n-2)!(2n-1)!}r_{n}^{\rm Bo}(b_n)
$$
\end{theorem} 

\begin{figure}[ht]
\centerline{\epsfbox{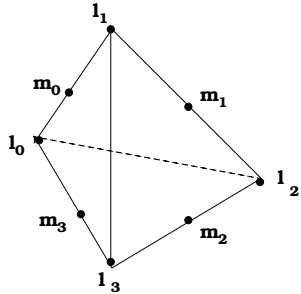}}
\caption{A special configuration of $8$ points in $\PP^{3}$.}
\label{fig1coc}
\end{figure}

{\bf 6. ${\Bbb P}^1 - \{0, \infty\}$ as a 
special stratum in the configuration space of $2n$ points in $\PP^{n-1}$ 
\cite{G4}}. 
A {\it special configuration}  is 
a configuration of $2n$  points  
\begin{equation} \label{--00}
(l_0,...,l_{n-1}, m_0,...,m_{n-1}) 
\end{equation} 
in $\PP^{n-1}$ such that  $l_0,...,l_{n-1}$ are vertices of a 
simplex in $\PP^{n-1}$ and $ m_i$ is a point 
on the edge $l_il_{i+1}$ of the simplex different from $l_i$ and 
$l_{i+1}$,  as on Figure \ref{fig1coc}. 

\begin{proposition} \label{6.12.02.1}
The set of special configurations of $2n$ points in $\PP^{n-1}$
is canonically identified with $\PP^1 \backslash \{0, \infty \}$. 
\end{proposition}

{\bf Construction}. Let $\widehat  m_i$
be the point of intersection of the line $l_il_{i+1}$ with the
hyperplane passing through all the points $m_j$ except $m_i$. 
Let $r(x_1, ..., x_4)$ be the cross-ratio of the four points on $\PP^1$. 
Let us define the generalized cross-ratio by \index{generalized cross-ratio} 
$$
r(l_0,...,l_{n-1},m_0,...,m_{n-1}) :=
r(l_i,l_{i+1},m_i,\widehat  m_{i+1}) \in F^*
$$ 
It does not depend on $i$, and provides  the desired isomorphism. 
Here is a different definition, which makes obvious 
the fact that the generalized cross-ratio is cyclically invariant. 
Consider the one dimensional subspaces $L_i, M_j$ in the $(n+1)$--dimensional 
vector space projecting to $l_i, m_j$. Then $L_i, M_i, L_{i+1}$ 
belong to a two dimensional subspace. The subspace $M_i$ 
provides a linear map $L_i \lra L_{i+1}$. The composition of these 
maps is a linear map $L_0 \lra L_0$. The element of $F^*$ describing this map   
is the generalized cross-ratio.

{\bf 7. Restriction of the Grassmannian $n$--logarithm to the special stratum}. 
The $n$--logarithm function $Li_n(z)$ has a  single-valued
version  (\cite{Z1}) 
\begin{eqnarray*} 
{\cal L}_{n}(z) &:=& \begin{array}{ll} 
{\rm Re} & (n:\ {\rm odd}) \\ 
{\rm Im} & (n: \ {\rm even}) \end{array} 
\left( \sum^{n-1}_{k=0} \beta_k
\log^{k}\vert z\vert \cdot Li_{n-k}(z)\right)\; , \quad n\geq 
2 \\ 
\end{eqnarray*}

It is continuous on $\C \PP^1$. Here   $\frac{2x}{e^{2x} -1}  =
\sum_{k=0}^{\infty}\beta_k x^k $, so $\beta_k = \frac{2^kB_k}{k!}$ 
where $B_k$ are Bernoulli numbers. For example ${\cal L}_2(z)$
is the Bloch - Wigner function.

Let us consider the following modification of the 
 function ${\cal L}_n(z)$ proposed by A. M. Levin in \cite{Le}:
\begin{equation} \label{6.28.04.1}
\widetilde {\cal L}_{n}(x):= 
\end{equation} 
$$
\frac{(2n-3)}{(2n-2)}
\sum_{\mbox{$k$ even;  $0 \leq k \leq n-2$}}
\frac{2^k (n-2)!(2n-k-3)!}{(2n-3)!(k+1)!(n-k-2)!} {\cal L}_{n-k}(x)\log^k|x|
$$

For example $\widetilde {\cal L}_{n}(x) =  {\cal L}_{n}(x)$ for $n \leq 3$, 
but already $\widetilde {\cal L}_{4}(x)$ is different from ${\cal L}_{4}(x)$. 
A direct integration carried out in Proposition 4.4.1 of \cite{Le} 
shows that  
$$
-(2\pi i)^{n-1}(-1)^{(n-1)(n-2)/2}\widetilde {\cal L}_{n}(x) =
$$
$$
\int_{\C \PP^{n-1}}\log|1-z_1| \prod_{i=1}^{n-1}d\log|z_i| \wedge \prod_{i=1}^{n-2}d\log|z_{i} - z_{i+1}| 
\wedge d\log|z_{n-1} -a| 
$$
This combined with Proposition \ref{0.3} below implies 
\begin {theorem}  \label{clpoly}
The value of function  ${\cal L}^G_n$ on special configuration
(\ref{--00}) equals 
$$
-(-1)^{n(n-1)/2}4^{n-1}{2n-2\choose n-1}^{-1}\widetilde {\cal L}_{n}(a);
\qquad a = r(l_0,...,l_{n-1},m_0,...,m_{n-1})
$$ 
\end {theorem}

{\bf 8. Computation of the Grassmannian $n$--logarithm}. 
It follows from Theorem \ref{clpoly} that ${\cal L}^{G}_2(l_1, ..., l_4)=-2 {\cal L}_2(r(l_1, ..., l_4))$.

It was proved in Theorem 1.3 of \cite{GZ} that 
$$
{\cal L}^G_3(l_0, ..., l_5) = \frac{1}{90} {\rm Alt}_6
{\cal L}_3 (r_3(l_0, ..., l_5))   
+
$$
\begin{equation} \label{3.31.03.5}
\frac{1}{9}{\rm Alt}_6\Bigl(\log 
|\Delta(l_0, l_1, l_2)| \log 
|\Delta(l_1, l_2, l_3)| \Delta(l_2, l_3, l_4)| \Bigr) 
\end{equation}
We will continue this discussion in chapter 5. 

The functions ${\cal L}^{G}_n$ for $n>3$ can not be expressed 
via classical polylogarithms.

\section{Polylogarithmic motivic complexes}

{\bf 1. The groups ${\cal B}_{n}(F)$ and 
polylogarithmic motivic complexes (\cite{G1}-\cite{G2})}.  
\index{Polylogarithmic motivic complexes} For  a set $X$ denote by  $\Z[X]$ 
the free abelian group generated by symbols $\{x\}$ where $x$ run
through all elements of the set $X$. Let $F$ be an arbitrary field. 
We  define inductively subgroups
${\cal R}_{n}(F)$ of $ \Z[P^{1}_{F}]$, $n\geq 1$ and set 
$$
{\cal B}_{n}(F):=\Z[P^{1}_{F}]/{\cal R}_{n}(F) 
$$ 
One has 
$$ 
{\cal R}_{1}(F):=(\{x\} + \{y\} - \{xy\},(x,y\in F^{\ast}); 
\{0\};\{\infty\})\;; \qquad   
{\cal B}_{1}(F) = F^{\ast}
$$ 
 Let $\{x\}_{n}$ be the image of $\{x\}$ in ${\cal 
B}_{n}(F)$. Consider homomorphisms
\begin{eqnarray} 
 \mbox{$\Z$}[P^{1}_{F}] 
\stackrel{\delta_{n}}{\longrightarrow}  \left\{ 
\begin{array}{lll} 
{\cal B}_{n-1}(F)\otimes F^{\ast} &:& n\geq 3 \\ 
\Lambda^{2}F^{\ast} &:& n=2\end{array}\right. 
\end{eqnarray}
\begin{eqnarray}
 \delta_{n}:\{x\}\mapsto \left\{ \begin{array}{lll} 
\{x\}_{n-1}\otimes x &:& n\geq 3 \\ 
(1-x)\wedge x &:& n=2\end{array}\right.
 \qquad \delta_{n}:\{\infty\},\{0\},\{1\}\mapsto 0 
\end{eqnarray} 
  Set 
$ {\cal A}_{n}(F):={\rm Ker}\ \delta_{n}\; . $ 
 Any element $\alpha(t) = \Sigma n_{i}\{f_{i}(t)\} \in \Z[P^{1}_{F(t)}]$ has a specialization $\alpha(t_{0}):=\Sigma 
n_{i}\{f_{i}(t_{0})\}\in \Z[P^{1}_{F}]$ at each point $t_{0}\in 
P^{1}_{F}$.  
 
\begin {definition}
\label {2.10}  
${\cal R}_{n}(F)$  is generated by 
elements $\{\infty\},\{0\}$  and   $\alpha(0)-\alpha(1)$  where $\alpha(t)$ runs 
through all elements of ${\cal A}_{n}(F(t))$.
\end {definition}
 Then $\delta_{n}\Bigl({\cal R}_{n}(F)\Bigl)=0$ (\cite{G1}, 1.16).  So we get 
homomorphisms  
$$
\delta_n: {\cal B}_{n}(F) \longrightarrow 
{\cal B}_{n-1}(F)\otimes F^{\ast}, \quad n\geq 3; \quad  \delta_2: {\cal
B}_{2}(F) \longrightarrow  \Lambda^{2}F^{\ast} 
$$
and finally the {\it polylogarithmic motivic complex} $\Gamma (F,n)$: 
$$
{\cal B}_{n}\stackrel{\delta}{\rightarrow} {\cal B}_{n- 
1}\otimes F^{\ast} \stackrel{\delta}{\rightarrow} {\cal B}_{n- 
2}\otimes \Lambda^{2}F^{\ast}\stackrel{\delta}{\rightarrow} \ldots \stackrel{\delta}{\rightarrow} 
{\cal B}_{2}\otimes \Lambda^{n- 
2}F^{\ast}\stackrel{\delta}{\rightarrow} \Lambda^{n}F^{\ast} 
$$ 
where   
$\delta:\{x\}_{p}\otimes \bigwedge^{n-p}_{i=1} y_{i}\to 
\delta_p(\{x\}_{p})\wedge \bigwedge^{n-p}_{i=1}y_{i} $ and ${\cal B}_{n}$ is in degree $1$. 

\begin{conjecture} \label{3.16.03.10}
$H^i\Gamma(F,n)\otimes \Q = {\rm gr}^{\gamma}_nK_{2n-i}(F)\otimes \Q$. 
\end{conjecture}
Denote by $\widehat {\cal L}_{n}$ the function 
${\cal L}_{n}$, multiplied by $i$ for even $n$ and unchanged for odd $n$. 
There is  a well defined homomorphism (\cite{G2}, Theorem 1.13):
$$
\widehat {\cal L}_{n}: {\cal B}_{n}(\C)) \lra \R(n-1); \qquad 
\widehat {\cal L}_{n}(\sum m_i \{z_i\}_n):= \sum m_i \widehat {\cal L}_{n}(z_i)
$$
There are canonical homomorphisms
\begin{equation} \label{4.15.03.1}
B_n(F) \lra {\cal B}_n(F); \quad \{x\}_n \lms \{x\}_n, \quad n=1,2,3.
\end{equation}
They are isomorphisms for $n=1,2$ and expect to be an isomorphism 
for $n=3$, at least modulo torsion. 
 
{\bf 2. The residue homomorphism for complexes $\Gamma(F,n)$ (1.14,  
\cite{G1})}. \index{The residue homomorphism} 
Let $F=K$ be a  field with 
a discrete valuation $v$,  the residue field  $k_v$ and  the group of units
$U$. Let $u \rightarrow \bar u$ be the projection $U \rightarrow
k_v^{\ast}$. Choose  a uniformizer $\pi$. There is a homomorphism $\theta:\Lambda^{n}K^{\ast} 
\longrightarrow\Lambda^{n-1} k_v^{\ast}$ uniquely defined 
by the  following properties $(u_{i}\in U)$: 
$$
\theta\; (\pi\wedge u_{1}\wedge \cdots\wedge u_{n-1}) 
= \bar u_{1}\wedge\cdots \wedge \bar u_{n-1}; \qquad  \theta\;
(u_{1}\wedge \cdots \wedge u_{n}) = 0 
$$
It is clearly independent  of $\pi$. Define a homomorphism $s_{v}:\Z[P^{1}_{K}]\longrightarrow \Z[P^{1}_{ k_v}]$ by setting
 $s_{v}\{ x\} =  \{ \bar x\}  \mbox{ if $x$ is a unit}$ and $0$  
 otherwise.
It induces a homomorphism 
$s_{v} :  {\cal B}_{m}(K)\longrightarrow {\cal B}_{m}(k_v)$. 
Put
$$ 
\partial_{v}:= s_{v}\otimes \theta : \quad {\cal B}_{m}(K)\otimes 
\Lambda^{n-m} K^{\ast}\quad \longrightarrow \quad {\cal B}_{m} (k_v) 
\otimes \Lambda^{n-m-1} k_v^{\ast} 
$$ 
It  defines a morphism 
of complexes
$\partial_{v} : \Gamma (K,n)\longrightarrow  \Gamma( k_v,
n-1)[-1]$.

{\bf 3. A variation of mixed  
Hodge structures on $P^1(\C) - \{0, 1, \infty\}$  
corresponding to the classical polylogarithm ${\rm Li}_n(z)$ (\cite{D2})}. 
Its fiber $H(z)$ over a point $z$ is  
described via the period matrix
\begin{equation} \label{4.20.01.1}
\left (\matrix{1&0&0&... &0 \cr 
{\rm Li}_1(z)& 2 \pi i &0 &... & 0 \cr 
{\rm Li}_2(z) & 2 \pi i \log z & (2 \pi i)^2& ... & 0\cr 
    ...  
& ...&    ... & ... &...\cr
{\rm Li}_n(z)& 2 \pi i \frac{\log^{n-1}z}{(n-1)!}&
(2 \pi i)^2 \frac{\log^{n-2}z}{(n-2)!} & ... & (2 \pi i)^n 
\cr}\right ) 
\end{equation}
Its entries are defined using   
analytic continuation to the point $z$ along a path $\gamma$ 
from a given point in $\C$ where all the entries are defined by
power series expansions, say the point $1/2$. 

Here is a more  natural way to define the entries. 
Consider the following regularized iterated integrals along a 
certain fixed path $\gamma$ between $0$ to $z$: 
$$
{\rm Li}_n(z) = \int_0^z\frac{dt}{1-t} \circ \underbrace{\frac{dt}{t}\circ ... \circ 
\frac{dt}{t}}_{n-1 \quad  \mbox{times}}; \qquad 
\frac{\log^{n}z}{n!} =   \int_0^z \underbrace{\frac{dt}{t} \circ  ... 
\circ\frac{dt}{t}}_{n \quad  \mbox{times}} 
$$
To regularize the divergent integrals we take the lower limit of integration to be $\varepsilon$.
Then it is easy to show that the integral has an asymptotic expansion of type $I_0(\varepsilon) + 
I_1(\varepsilon)\log  \varepsilon + ... + I_k(\varepsilon)\log^k  \varepsilon$,  
where all the functions  $I_i(\varepsilon)$ are smooth at $\varepsilon =0$. 
Then we take $I_0(0)$ to be the regularized value.

Now let us define a mixed Hodge structure  $H(z)$. 
Let $\C^{n+1}$ be the standard vector space with basis 
$(e_0, ..., e_n)$, and $V_{n+1}$ the $\Q$-vector subspace 
spanned by the columns of the matrix (\ref{4.20.01.1}). 
Let $W_{-2(n+1-k)}V_{n+1}$ be the subspace spanned by  the first $k$ columns, 
counted from the right to the left. One shows that these subspaces do not
depend on the choice of the path $\gamma$, i.e. they are 
well-defined inspite of the 
multivalued nature of the entries of the matrix. Then $W_{\bullet}V_{n+1}$ is
 the weight filtration. 
We define the  Hodge filtration $F^{\bullet}\C^{n+1}$ by setting $F^{-k}\C^{n+1}:= \langle e_0, ..., e_k\rangle_{\C}$. 
It is opposite to the weight filtration. We get a Hodge-Tate structure, i.e. $h^{pq}=0$ unless $p \not = q$. One checks that the family of Hodge-Tate structures $H(z)$ forms a unipotent variation of 
Hodge-Tate structures on $P^1(\C) - \{0, 1, \infty\}$. 

Let $n \geq 0$. An $n$-framed  
Hodge-Tate structure  $H$ is a triple $(H, v_0, f_n)$, where 
$v_{0} : \Q(0) \longrightarrow 
gr^W_{0}H$ and   $f_{n}: gr^W_{-2n} H \longrightarrow \Q(n)$ are nonzero morphisms. 
A framing plus a choice of a splitting of the weight filtration determines a period of a 
Hodge-Tate structure.  
  Consider the  coarsest equivalence
relation on the set of all $n$-framed Hodge-Tate structures  for which $M_1 \sim M_2$
if there is a map $ M_1 \to  M_2$ respecting the frames. Tnen the set ${\cal H}_n$ of the equivalence classes has a natural abelian group structure.  Moreover
\begin{equation} \label{HH}
{\cal H}_{\bullet} := \oplus_{n \geq 0}{\cal H}_n
\end{equation}
has a natural Hopf algebra structure with a coproduct $\Delta$, see the Appendix of \cite{G12}. 

Observe that ${\rm Gr}^W_{-2k}H_n(z) = \Q(k)$    
for  $-n \leq k \leq 0$. 
Therefore $H_n(z)$ has a natural framing such that the corresponding period 
is given by the function ${\rm Li}_n(z)$. The obtained framed object is denoted 
by ${\rm Li}^{\cal H}_n(z)$. To define it for $z = 0, 1, \infty$ we use the 
specialization functor to the punctured tangent space at $z$, and then take the fiber 
over the tangent vector corresponding to the parameter $z$ on $P^1$. 
It is straitforward to see that the coproduct $\Delta {\rm Li}^{\cal H}_n(z)$ is computed
by the formula
\begin{equation} \label{cani5d}
\Delta {\rm Li}^{\cal H}_n(z)  = \sum_{k=0}^n {\rm Li}^{\cal H}_{n-k}(z) \otimes 
\frac{\log^{\cal H}(z)^k}{k!}
\end{equation}
where $\log^{\cal H}(z)$ is the $1$-framed Hodge-Tate structure corresponding to $\log(z)$. 

{\bf 4. A motivic proof of the weak version of Zagier's conjecture}. Our goal is the following result, which was proved in \cite{dJ3} and in the unfinished manuscript \cite{BD2}. 
The proof below uses a different 
set of ideas. It follows the framework described in Chapter 13 of \cite{G8}, and quite 
close to the approach outlined in \cite{BD1}, although it is formulated a bit differently, 
using the polylogarithmic motivic complexes.  

\begin{theorem} \label{6.28.04.10} Let $F$ be a number field. Then there exists 
a homomorphism
$$
l_n: H^1(\Gamma(F, n)\otimes \Q) \lra K_{2n-1}(F) \otimes \Q
$$
such that for any embedding $\sigma: F \hookrightarrow  \C$ one has the following commutative diagram
\begin{equation} \label{6.29.04.10}
\begin{array}{ccccc}
H^1(\Gamma(F, n)\otimes \Q)&\stackrel{\sigma_* }{\lra} & H^1(\Gamma(\C, n)\otimes \Q)&
\stackrel{{\cal L}_n }{\lra} &\R(n-1)\\
&&&&\\
\downarrow l_n&& && \downarrow =\\
&&&&\\
K_{2n-1}(F) \otimes \Q& \stackrel{\sigma_* }{\lra}& K_{2n-1}(\C) \otimes \Q &
\stackrel{r_{Bo}}{\lra}&\R(n-1)
\end{array}
\end{equation}
\end{theorem}

{\bf Proof}. We will use the following background material and facts:

i) The existence of the abelian tensor category ${\cal M}_T(F)$ of mixed Tate motives 
over a number field $F$, satisfying all the desired properties, including 
Beilinson's formula expressing the Ext groups via the rational $K$-theory of $F$ 
and 
the Hodge realization functor. 
See \cite{DG} and the references there.

ii) The formalism of 
mixed Tate categories, including the description of the fundamental Hopf algebra 
${\cal A}_{\bullet}({\cal M})$ 
of a mixed Tate category via  framed objects in ${\cal M}$. See \cite{G12}, Section 8. 
The fundamental Hopf algebra of the category ${\cal M}_T(F)$ is denoted 
${\cal A}_{\bullet}(F)$. For example for the category of mixed Hodge-Tate structures 
fundamental Hopf algebra is the one ${\cal H}_{\bullet}$ from (\ref{HH}). 
Let 
$$
\Delta: {\cal A}_{\bullet}(F) \lra {\cal A}_{\bullet}(F)^{\otimes 2}
$$
be the coproduct in the Hopf algebra ${\cal A}_{\bullet}(F)$, 
and $\Delta': =  \Delta - {\rm Id}\otimes 1 + 1 \otimes {\rm Id}$ is the restricted 
coproduct. 
The key fact is a canonical isomorphism
\begin{equation} \label{cani}
{\rm Ker}\Delta' \cap {\cal A}_{n}(F) \stackrel{\sim}{=} K_{2n-1}(F)\otimes \Q
\end{equation} 
Since ${\cal A}_{\bullet}(F)$ is graded 
by $\geq 0$ integers, and ${\cal A}_0(F) = \Q$, formula  (\ref{cani}) for $n=1$
 reduces to an isomorphism 
\begin{equation} \label{cani1}
{\cal A}_1(F) \stackrel{\sim}{=} F^* \otimes \Q
\end{equation}

iii) The existence of the motivic classical polylogarithms
\begin{equation} \label{cani2}
{\rm Li}^{\cal M}_n(z) \in {\cal A}_n(F), \qquad z \in P^1(F)
\end{equation}
They were defined in Section 3.6 of \cite{G12} using either a geometric construction 
of \cite{G9}, or a construction of the motivic fundamental torsor of 
path between the tangential base points given in \cite{DG}.  
In particular one has ${\rm Li}^{\cal M}_n(0)  = {\rm Li}^{\cal M}_n(\infty)
=0$. 
A natural construction of the elements (\ref{cani2}) using the moduli space $\overline
{\cal
  M}_{0, n+3}$ is given in Chapter 4.6 below.

iv) The Hodge realization of the element (\ref{cani2}) 
is equivalent to the framed Hodge-Tate structure ${\rm Li}^{\cal H}_n(z)$ 
from  Chapter 4.3. This fact is more or less straitforward 
if one uses the fundamental torsor of path on the 
punctured projective line to define the element ${\rm Li}^{\cal M}_n(z)$, and 
follows from the general specialization theorem proved in 
\cite{G9} if one uses the approach of loc. cit. 
This implies that the Lie-period 
of ${\rm Li}^{\cal H}_n(z)$ equals to ${\cal L}_n(z)$. Indeed, for the 
Hodge-Tate structure $H(z)$ assigned to ${\rm Li}_n(z)$ this was shown in \cite{BD1}. 

v) The crucial formula (Section 6.3 of \cite{G12}):
\begin{equation} \label{cani5}
\Delta {\rm Li}^{\cal M}_n(z)  = \sum_{k=0}^n {\rm Li}^{\cal M}_{n-k}(z) \otimes 
\frac{\log^{\cal M}(z)^k}{k!}
\end{equation}
where $\log^{\cal M} (z) \in {\cal A}_1(F)$ is the element 
corresponding to $z$ under the isomorphism (\ref{cani}), and 
$\log^{\cal M}(z)^k \in {\cal A}_k(F)$ is its 
$k$-th power. It follows from formula (\ref{cani5d}) using the standard trick based on Borel's theorem to reduce a motivic claim to the corresponding Hodge one. 

vi) The Borel regulator map on $K_{2n-1}(F)\otimes \Q$, which sits via (\ref{cani}) inside of 
${\cal A}_{n}(F)$, is induced by the Hodge realization functor on 
the category ${\cal M}_T(F)$. 

Having this background, we proceed as follows. Namely, let 
$$
{\cal L}_{\bullet}(F) := \frac{{\cal A}_{>0}(F)}{{\cal A}_{>0}(F)^2}
$$
be the fundamental Lie coalgebra of ${\cal M}_T(F)$. Its cobracket $\delta$ is induced by $\Delta$. Projecting the element 
(\ref{cani2}) into ${\cal L}_{n}(F)$ we get an element $l^{\cal M}_n(z) \in {\cal L}_n(F)$ such that 
\begin{equation} \label{cani3}
\delta l^{\cal M}_n(z) = l^{\cal M}_{n-1}(z) \wedge z 
\end{equation}

Consider the following map 
$$
\widetilde l_n: \Q[P^1(F)] \lra {\cal L}_n(F), \qquad \{z\} \lms l^{\cal M}_n(z)
$$
\begin{proposition} \label{6.20.4.11}
The map $\widetilde l_n$ kills the subspace ${\cal R}_n(F)$, 
providing a well defined homomorphism
$$
l_n: {\cal B}_n(F)  \lra {\cal L}_n(F) ; \qquad \{z\}_n \lra \widetilde l_n(z)
$$
\end{proposition}

{\bf Proof}. 
We proceed by the induction on $n$. The case $n=1$ is self-obvious. 
Suppose we are done for $n-1$. 
Then there is the  following commutative diagram:  
$$
\begin{array}{ccc}
\Q[F] & \stackrel{\delta_n}{\lra} &{\cal B}_{n-1}(F)\otimes F^*\\
&&\\
\downarrow \widetilde l_n &&\downarrow l_{n-1}\wedge {\rm Id} \\
&&\\
{\cal L}_n(F) & \stackrel{\delta}{\lra} &\oplus_{k\leq n/2}{\cal L}_{n-k}(F)\wedge {\cal L}_{k}(F)\
\end{array}
$$
 Indeed, its commutativity is equivalent to the basic formula (\ref{cani3}).

Let $x \in P^1(F)$. Recall the specialization at $x$ homomorphism
$$
s_x: {\cal B}_n(F(T)) \lra {\cal B}_n(F), \qquad \{f(T)\}_n \lms \{f(x)\}_n
$$
It gives rise to the specialization homomorphism 
$$
s_x': {\cal B}_{n-1}(F(T))\otimes F(T)^* \lra {\cal B}_{n-1}(F)\otimes F^*, 
$$
$$
\{f(T)\}_{n-1}\otimes g(T) \lms 
\{f(x)\}_{n-1}\otimes \frac{g(T)}{(T-x)^{v_x(g)}}(x)
$$
(Use the local parameter $T^{-1}$ when $x = \infty$).

Now let 
$$
\alpha(T) \in {\rm Ker}\Bigl(\delta_n: \Q[F(T)] \to {\cal B}_{n-1}(F(T))\otimes F(T)^*\Bigr)
$$
Using ${\rm Li}^{\cal M}_n(0) = {\rm Li}^{\cal M}_n(\infty) = 0$, 
for any $x \in P^1(F)$ one has $\delta(\widetilde l_n(\alpha(x))) =0$. Thus 
$$
\widetilde l_n(\alpha(x)) \in K_{2n-1}(F)\otimes \Q \subset {\cal L}_n(F)
$$
Let us show that this element is zero. Given an embedding $\sigma: F \hookrightarrow
 \C$, 
write $\sigma(\alpha(T)) = 
\sum_i n_i\{f^{\sigma}_i(T)\}$.  
Applying the 
Lie-period map to this element and using v) we get  $\sum_i n_i{\cal L}_nf^{\sigma}_i(z)$. 
By Theorem 1.13 in \cite{G2} the condition on $\alpha(T)$ implies that 
this function is constant on $\CP^1$. 
Thus the difference of its values at $\sigma(x_1)$ and $\sigma(x_2)$, where 
$x_0, x_1 \in P^1(F)$, is zero. On the other hand thanks to v) and vi) it coincides with the Borel regulator map applied to the corresponding element of $K_{2n-1}(F)\otimes \Q$. 
Thus  the injectivity of the Borel 
regulator map proves the claim. 
So 
$\widetilde l_n(\alpha(x_0) - \alpha(x_1)) =0$. The proposition is proved. 

Proposition \ref{6.20.4.11} implies that we get a homomorphism of complexes
$$
\begin{array}{ccc}
{\cal B}_n(F) & \stackrel{\delta}{\lra} &{\cal B}_{n-1}(F)\otimes F^*\\
&&\\
\downarrow l_n &&\downarrow l_{n-1}\wedge {\rm Id} \\
&&\\
{\cal L}_n(F) & \stackrel{\delta}{\lra} &\oplus_{k\leq n/2}{\cal L}_{n-k}(F)\wedge {\cal L}_{k}(F)\
\end{array}
$$

The theorem follows immediately from this. Indeed, it remains to 
check commutativity of the diagram (\ref{6.29.04.10}), and it follows from v) and vi). 


If we assume the existence of the hypothetical abelian category of 
mixed Tate motives over an arbitrary field $F$, the same argumentation as above 
(see Chapter 6.1) 
implies the following result: one should have canonical 
homomorphisms 
$$
H^i(\Gamma(F, n))\otimes \Q \lra {\rm gr}^{\gamma}_nK_{2n-i}(F)\otimes \Q
$$
The most difficult part of Conjecture \ref{3.16.03.10} 
says that this maps are supposed to be isomorphisms. 

In the next section we define a regulator map on the polylogarithmic motivic complexes. 
Combined with these maps, it should give an explicit construction of the regulator map.

{\bf 5. A construction of the motivic $\zeta$-elements (\ref{2.11.04.1})}. 
The formula (\ref{12q}) leads to the motivic extension $\zeta^{\cal M}(n)$ as
follows (\cite{GoM}). Recall the moduli space $\overline {\cal M}_{n+3}$ 
parametrising stable curves of genus zero with $n+3$ marked points. 
It contains as an open subset the space ${\cal M}_{n+3}$ parametrising 
the $(n+3)$-tuples of distinct points on ${\Bbb P}^1$ modulo ${\rm Aut}({\Bbb
  P}^1)$. Then the complement $\partial {\cal M}_{n+3}:= \overline {\cal M}_{n+3} - {\cal M}_{n+3}$
is a normal crossing divisor, and the pair  $({\cal M}_{n+3}, \partial {\cal
  M}_{n+3})$ 
is defined over $\Z$. Let us identify sequences $(t_1, ..., t_n)$ of distinct
complex numbers  different from $0$ and $1$  
with the points $(0, t_1, ..., t_n, 1, \infty)$ of ${\cal M}_{n+3}(\C)$. 
Let us consider the integrand in (\ref{12q}) as a holomorphic form 
on ${\cal M}_{n+3}(\C)$. Meromorphically extending it to 
 $\overline {\cal M}_{n+3}$ we get a differential form with logarithmic
 singularities $\Omega_n$. Let $A_n$ be its divisor.  
Similarly, embed the integration simplex ${0 < t_1 < ... < t_n < 1}$ into the
set of real points of $\overline {\cal M}_{n+3}(\R)$, take its closure
$\Delta_n$ 
there, and consider  the Zariski closure $B_n$ of its boundary $\partial \Delta_n$. Then one can
show that the mixed motive
\begin{equation} \label{3435}
H^n(\overline {\cal M}_{n+3} - A_n, B_n - (A_n\cap B_n))
\end{equation}
is a mixed Tate motive over ${\rm Spec}(\Z)$. Indeed, it is easy to prove that 
its $l$-adic realization 
is unramified outside $l$, and is glued from the Tate modules of different
weights, and then refer to \cite{DG}.  The mixed motive (\ref{3435}) comes equipped with an
additional data, framing, given by non-zero morphisms
\begin{equation} \label{7.14.04.1}
[\Omega_n]: \Z(-n) \lra {\rm gr}^W_{2n}H^n(\overline {\cal M}_{n+3} - A_n, B_n
- (A_n\cap B_n)), 
\end{equation}
\begin{equation} \label{7.14.04.2}
[\Delta_n]:  {\rm gr}^W_{0}H^n(\overline {\cal M}_{n+3} - A_n, B_n
- (A_n\cap B_n))\lra \Z(0), 
\end{equation}
There exists the minimal subquotient of the mixed motive (\ref{3435}) which
inherits non-zero framing. It delivers the extension class $\zeta^{\cal M}(n)$.   
Leibniz formula 
(\ref{12q}) just means that $\zeta(n)$ is its period.

\begin{figure}[ht]
\centerline{\epsfbox{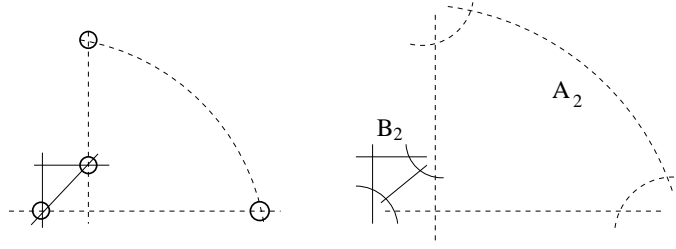}}
\caption{Constructing $\zeta^{\cal M}(2)$}
\label{K1}
\end{figure}

{\bf Example}. To construct $\zeta^{\cal M}(2)$, take the pair of triangles in $P^2$ 
shown on the left of Figure \ref{K1}. The triangle shown by the punctured lines is the 
divisor of poles of the differential $d\log(1-t_1) \wedge d\log t_2$ 
in (\ref{12q}), and the second triangle is 
the algebraic closure of the boundary of the integration cycle $0 \leq t_1 \leq t_2 \leq 1$. 
 The corresponding configuration of six lines is defined uniquely 
up to a projective equivalence. Blowing up the four points shown by little circles on 
Figure \ref{K1} (they are the triple intersection points of 
the lines), 
we get the moduli 
space $\overline {\cal M}_{0,5}$. Its boundary is the union of the two 
pentagons, $A_2$ and $B_2$, projecting to the two 
triangles on $P^2$. Then $\zeta^{\cal M}(2)$ is 
a subquotient of the mixed Tate motive 
$H^2(\overline {\cal M}_{0,5} - A_2, B_2 - (A_2 \cap B_2))$ over ${\rm Spec}(\Z)$.   
Observe that an attempt to use a similar construction for 
the pair of triangles in $P^2$ fails since 
there is no non-zero morphism $[\Delta_2]$ in this case. 
Indeed, 
there are  two vertices of the
$B$-triangle shown on the left of  Figure \ref{K1} lying at the sides of
the $A$-triangle, and therefore the chain $0 \leq t_1 \leq t_2 \leq 1$ does not give rise 
to a relative class in $H_2(P^2(\C) - A, B - (A\cap B))$. 

{\bf 6. A geometric construction of the motivic classical polylogarithm ${\rm
    Li}^{\cal M}_n(z)$}. 
The above construction is easily generalized to the case of the classical
    polylogarithm. Let $A_n(z)$ be the divisor of the meromorphic differential
    form 
$$
\Omega_n(z) := \frac{dt_1}{z^{-1}- t_1}\wedge \frac{dt_2}{t_2}\wedge ... \wedge \frac{dt_n}{t_n}
$$
extended as a rational form to $\overline {\cal M}_{0, n+3}$. 
Suppose that $F$ is a number field and $z \in F$. Then there is  the following
mixed Tate motive over $F$:  
\begin{equation} \label{3435s}
H^n(\overline {\cal M}_{n+3} - A_n(z), B_n - (A_n(z)\cap B_n))
\end{equation}
One 
checks that  the vertices (that is the
zero-dimensional strata) of the divisor $B_n$ are disjoint with the divisor 
$A_n(z)$.  
Using this we define a framing $([\Omega_n(z)], [\Delta_n])$ 
on the mixed motive (\ref{3435s}) 
similar to the one (\ref{7.14.04.1}) -
(\ref{7.14.04.2}). The geometric condition on the divisor  $A_n(z)$ is used 
to show that the framing morphism $[\Delta_n]$ is non-zero. We define ${\rm
    Li}^{\cal M}_n(z)$ as the framed mixed Tate motive (\ref{3435s}) with the  
  framing $([\Omega_n(z)], [\Delta_n])$.  
A similar construction for the multiple polylogarithms was worked out in the
Ph. D. thesis of Q. Wang \cite{Wa}.

It follows from a general result in \cite{G9} that the Hodge realization of 
${\rm
    Li}^{\cal M}_n(z)$ is equivalent to the framed mixed Hodge structure 
${\rm Li}^{\cal H}_n(z)$. 
If $n=2$ the two corresponding mixed Hodge structures are isomorphic 
(\cite{Wa}), while in general they are not isomorphic but equivalent 
as framed  mixed Hodge structures.

\section{Regulator maps on the polylogarithmic motivic complexes}

In this chapter we define explicitly 
these regulator maps via classical polylogarithm functions 
following \cite{G6} and \cite{G2}. This implies how the special 
values of $\zeta$-functions of algebraic varieties outside of the 
critical strip  should be expressed using the classical polylogarithms. 

{\bf 1. The numbers $\beta_{k,p}$}. 
Define for any integers $p \geq 1$ and $k \geq 0$ the numbers 
$$
\beta_{k, p}:= \quad (-1)^p (p-1)! \sum_{0 \leq i \leq [\frac{p-1}{2}]} \frac{1}{(2i+1)!} \beta_{k+p-2i}
$$
For instance
$\beta_{k, 1} = -\beta_{k+ 1}; \quad \beta_{k, 2} = \beta_{k+ 2}; \quad \beta_{k, 3} = 
-2\beta_{k+ 3} - \frac{1}{3}\beta_{k+ 1}$. 

One has  recursions
\begin{equation} \label{12.99}
2p\cdot \beta_{k+1, 2p} =  - \beta_{k, 2p+1} -\frac{1}{2p+1}\beta_{k+1}; \qquad 
(2p - 1)\cdot \beta_{k+1, 2p-1} =  - \beta_{k, 2p} 
\end{equation}
These recursions together with  $\beta_{k, 1} = -\beta_{k+ 1}$ determine the numbers $\beta_{k, p}$. 

Let $m \geq 1$. Then one can show that 
$$ 
\beta_{0, 2m} = \beta_{0, 2m+1} =  \frac{1}{2m+1}, \quad  
\beta_{1, 2m-1}  =  -\frac{1}{(2m-1)(2m+1)}, \quad \beta_{1, 2m}  = 0
$$

{\bf 2. The regulator map on the polylogarithmic motivic complexes 
in the case $F= \C(X)$, where $X$ is  a complex  algebraic variety $X$}. 
Let us define differential $1$-forms $\widehat {\cal L}_{p,q}$ on $\CP^1 \backslash \{0, 1, \infty\}$ for $q\geq 1$  as follows:
\begin{equation} \label{TG10}
\widehat {\cal L}_{p,q}(z):= \widehat {\cal L}_{p}(z) \log^{q-1}|z| \cdot d\log |z|, \quad p \geq 2
\end{equation}
$$
\widehat {\cal L}_{1,q}(z):=  \alpha(1-z, z) \log^{q-1}|z|
$$
For any rational function $f$ on a complex variety $X$ 
the $1$-form $\widehat {\cal L}_{p,q}(f)$ provides  a distribution on $X(\C)$. 
Set 
$$
{\cal A}_m\left\{\bigwedge_{i=1}^{2p}d\log|g_i| \wedge \bigwedge_{i=2p+1}^{m}di \arg g_{j}
  \right\}:= 
$$
$$
{\rm Alt}_m\left\{\frac{1}{(2p)!(m-2p)!}\bigwedge_{i=1}^{2p}d\log|g_i| \wedge \bigwedge_{i=2p+1}^{m}di \arg g_{j} \right\} 
$$
and 
$$
{\cal A}_m\left\{\log|g_1| \cdot \bigwedge_{i=2}^{p}d\log|g_i| \wedge \bigwedge_{i=p+1}^{m}di \arg g_{j} \right\}:= 
$$
$$
{\rm Alt}_m\left\{\frac{1}{(p-1)!(m-p)!}\log|g_1|  \cdot \bigwedge_{i=2}^{p}d\log|g_i| \wedge \bigwedge_{i=p+1}^{m}di \arg g_{j}   \right\}
$$
So
${\cal A}_m(F(g_1, .., g_m))$ is a weighted alternation 
(we divide by the order of the stabilizer of the term we alternate).

Let $f, g_1, ..., g_m$ be rational functions on a complex variety $X$. 
Set
$$
r_{n+m}(m+1):\{f\}_n \otimes g_1 \wedge ... \wedge g_m \lms 
$$
\begin{equation} \label{12.26.11}
\widehat {\cal L}_{n}(f) \cdot {\cal A}_m\left\{  \sum_{p \geq 0 } \frac{1}{2p+1}
 \bigwedge_{i=1}^{2p} d\log|g_i| \wedge \bigwedge_{j=2p+1}^{m} d i \arg g_{j}
 \right\} +
\end{equation}
\begin{equation} \label{12.26.12}
\sum_{k \geq 1 }\sum_{1 \leq p \leq m} \beta_{k, p} \widehat {\cal L}_{n-k,k}(f) \wedge 
{\cal A}_m\left\{\log|g_1| \bigwedge_{i=2}^p d \log|g_{i}| \wedge \bigwedge_{j= p+1}^m 
d i \arg g_{j}
 \right\} 
\end{equation}

\begin{proposition} \label{3.16.03.1}
The differential form $r_{n+m}(m+1)
(\{f\}_n \otimes g_1 \wedge ... \wedge g_m)$ defines a distribution on $X(\C)$.
\end{proposition}

{\bf Example 1}. $r_n(1)(\{f\}_n) = \widehat{{\cal L}}_n(f)$.

{\bf Example 2}. $r_n(n)(g_1 \wedge ... \wedge g_n) = 
r_{n-1}(g_1 \wedge ... \wedge g_n)$.

{\bf Example 3}. $m=1$,  $n$ is arbitrary. Then 
$$
r_{n+1}(2): \{f\}_n \otimes g \lms \quad \widehat {\cal L}_{n}(f) d i \arg g
- \sum_{k=1}^{n-1} \beta_{k+1} \widehat {\cal L}_{n-k, k}(f) \cdot \log|g| 
$$

 {\bf Example 4}. $m=2$, $n$ is arbitrary. 
$$
r_{n+2}(3): \{f\}_n \otimes g_1 \wedge g_2 \lms \widehat {\cal L}_{n}(f) \left\{  d i \arg g_1 \wedge d i \arg g_2 + 
\frac{1}{3} d \log |g_1| \wedge d \log |g_2|
\right\}
$$
$$
- \sum_{k=1}^{n-1} \beta_{k+1} \widehat {\cal L}_{n-k, k}(f) \wedge (\log|g_1| d i\arg g_2 -\log|g_2| d i\arg g_1)  
$$
$$
+ \sum_{k \geq 1} \beta_{k+2} \widehat {\cal L}_{n-k, k}(f) \wedge (\log|g_1| d \log|g_2| -\log|g_2| d \log|g_1|) 
$$

Let ${\cal A}^i(\eta_X)$ be the space of real 
smooth $i$-forms at the generic point $\eta_X:= {\rm Spec} \C(X)$ of a 
complex variety $X$. Let 
${\cal D}$ be  
the de Rham differential on distributions on $X(\C)$, 
and $d$ the de Rham differential on 
${\cal A}^i(\eta_X)$. For example:
$$
d \Bigl( d i \arg z\Bigl) = 0; \quad {\cal D}\Bigl( d i \arg z\Bigl) = 2 \pi i \delta(z)
$$
 
Recall the residue homomorphisms defined in Chapter 4.2. 
\begin{theorem} \label{CANHH}
\label {1.11a} a) The maps $r_n(\cdot)$ 
provide a homomorphism of complexes
$$
\begin{array}{ccccccc}
{\cal B}_{n}(\C(X))&\stackrel{\delta}{\rightarrow}& {\cal B}_{n- 
1}(\C(X))\otimes \C(X)^{\ast} & \stackrel{\delta}{\rightarrow} & \ldots &\stackrel{\delta}{\rightarrow} &
 \bigwedge^{n}\C(X)^{\ast}\\ 
&&&&&&\\
\downarrow r_n(1)& &\downarrow r_n(2)& & & &\downarrow r_n(n)\\
&&&&&&\\
{\cal A}^0(\eta_X)(n-1) &\stackrel{d}{\rightarrow}& {\cal A}^1(\eta_X)(n-1) &\stackrel{d}{\rightarrow}& ...
&\stackrel{d}{\rightarrow}& {\cal A}^{n-1}(\eta_X)(n-1)\\
\end{array}
$$

b) The maps $r_n(m)$ are compatible with the residues:
$$
{\cal D} \circ r_n(m) - r_n(m+1) \circ \delta = 2 \pi i \cdot \sum_{Y \in X^{(1)}}r_{n-1}(m-1) \circ \partial_{v_Y}, \quad m<n
$$
$$
{\cal D} \circ r_n(n) - \pi_n(d\log f_1 \wedge
... \wedge d\log f_n) = 
2 \pi i\cdot \sum_{Y \in X^{(1)}}r_{n-1}(n-1) \circ \partial_{v_Y}
$$
where $v_Y$ is  the  valuation  
on the field $\C(X)$ defined by a divisor $Y$.
 \end {theorem}

Let $X$ be a regular variety over $\C$.  Recall the $n$-th Beilinson-Deligne complex 
$\underline {\R}_{{\cal D}}(n)_X$ 
defined as a total complex associated with the following
bicomplex of sheaves in classical topology on $X(\C)$: 
$$
\begin{array}{ccccccccccc} \label{del}
\Bigl(\underline {\cal D}_{X}^{0}&\stackrel{d}{\longrightarrow}&\underline {\cal
D}_{X}^{1}&\stackrel{d}{\longrightarrow}&\ldots&\stackrel{d}{\longrightarrow}&\underline {\cal
D}^{n}_{X}&\stackrel{d}{\longrightarrow}&\underline {\cal
D}_{X}^{n+1}&\stackrel{d}{\longrightarrow}&\ldots\Bigr) \otimes \R(n-1)\\
&&&&&&\uparrow\pi_{n} &&\uparrow\pi_{n}&&\\
&&&&&&\Omega^{n}_{X, \log}
&\stackrel{\partial}{\longrightarrow}&\Omega_{X, \log}^{n+1}&\stackrel{\partial}{\longrightarrow}&
\end{array}
$$
Here $\underline{\cal D}^{0}_{X}$ is in
degree 1 and  
$(\Omega^{\bullet}_{X, \log}, \partial)$ is the de Rham complex of 
holomorphic forms with logarithmic singularities
 at infinity. We will denote by $\underline 
{\R}_{{\cal D}}(n)(U)$ the complex of the global sections. 

Theorem \ref{CANHH} can be reformulated as follows. 
Set $\widetilde r_n(i):= r_n(i)$ for $i<n$ and
$$
\widetilde r_n(n): \Lambda^n\C(X)^* \lra {\cal A}^{n-1}(\eta_X)(n-1) \oplus 
\Omega^n_{\log}(\eta_X)
$$
\begin{equation} \label{3/6/00.01}
f_1 \wedge ... \wedge f_n \lms r_n(n)(f_1 \wedge ... \wedge f_n) + d\log f_1 \wedge ... \wedge 
d\log f_n
\end{equation}

\begin{theorem} \label{CANHH*}Let $X$ be a complex algebraic variety. 
Then there is a homomorphism of complexes
\begin{equation} \label{CANHH**}
\widetilde r_n(\cdot): \Gamma(\C(X); n) \lra \underline 
{\R}_{{\cal D}}(n)({\rm Spec}\C(X))
\end{equation}
compatible with the residues 
as explained in the part b) of Theorem \ref{CANHH}. 
 \end {theorem}

{\bf 3. The general case}. Let $X$ be a regular projective variety over a
field $F$. 
Let $d:= {\rm dim} X$. Then the complex $\Gamma(X;n)$ {\it should be} defined as  the total complex 
of the following bicomplex:
$$
\Gamma(F(X); n) \lra \oplus_{Y_1 \in X^{(1)}}\Gamma(F(Y_1); n-1)[-1] \lra 
$$
$$
\oplus_{Y_2 \in X^{(2)}}\Gamma(F(Y_2); n-2)[-1] \lra ... 
\lra 
\oplus_{Y_d \in X^{(d)}}\Gamma(F(Y_d); n-d)[-d]
$$
where the arrows are provided by the residue maps, see \cite{G1}, p 239-240. 
The complex $\Gamma(X;n)\otimes \Q$ should be 
quasiisomorphic to the weight $n$ motivic complex. 

However there is a serious difficulty in the 
 definition of  the complex $\Gamma(X;n)$ 
for a general variety $X$ and $n>3$, (\cite{G1}, p. 240 ). 
It would be resolved if homotopy invariance 
of the polylogarithmic complexes  were known (Conjecture 1.39 in \cite{G1}). 
As a result we have an unconditional definition of the polylogarithmic 
complexes $\Gamma(X;n)$ only in the following cases: 

a) $X = Spec(F)$, $F$ is an arbitrary field.

b) $X$ is an regular curve over any field, and $n$ is arbitrary.

c) $X$ is an arbitrary regular scheme, but $n \leq 3$. 

Now let $F$ be a subfield of $\C$. Having in mind applications to arithmetic, 
 we will restrict ourself by the
case when $F = \Q$. 
Assuming we are working with one of the above cases, 
or assuming the above difficulty has been resolved, let us 
define the regulator map 
$$
\Gamma(X;n) \lra {\cal C}_{{\cal D}}(X(\C); \R(n))
$$
We specify it for each of 
$\Gamma(\Q(Y_k); n-k)[-k]$ where $k = 0, ..., d$. 
Namely, we take the homomorphism $r_{n-k}(\cdot)$ 
for ${\rm Spec}(\Q(Y_k)$ and multiply it by 
$(2 \pi i)^{n-k} \delta_{Y_{n-k}}$. Notice that the distribution $\delta_{Y}$ 
depends only on the generic point of a subvariety $Y$. 
Then   Theorem \ref{CANHH}, and in particular its part b), providing 
{\it compatibility with the residues} property, guarantee  
that we get a homomorphism of complexes. 
Here are some examples.

{\bf 4. The weight one}. The regulator map on the weight one motivic complex 
looks as follows: 
$$
\begin{array}{ccc}
\Q(X)^*& \lra &\oplus_{Y \in X^{(1)}}\Z\\
&&\\
\downarrow r_1(1)&&\downarrow r_1(2)\\
&&\\
{\cal D}^{0,0}&
\stackrel{2 \overline \partial \partial}{\lra} &{\cal D}_{cl, \R}^{1,1}(1)
\end{array}
$$
$$
r_1(2): Y_1 \lms 2 \pi i \cdot \delta_{Y_1}, \qquad 
r_1(1): f \lms \log|f|
$$
Here the top line is the weight $1$ 
motivic complex, sitting in degrees $[1,2]$.

{\bf 5. The weight two}. The regulator map on the weight two motivic complexes looks as follows. 
$$
\begin{array}{ccccccc}
{\cal B}_2(\Q(X))& \stackrel{\delta}{\lra} &\Lambda^2\Q(X)^*&\stackrel{\partial}{\lra} &\oplus_{Y \in X^{(1)}}\Q(Y)^*&\stackrel{\partial}{\lra}&\oplus_{Y \in X^{(2)}}\Z\\
&&&&&&\\
\downarrow r_2(1)&&\downarrow r_2(2)&&\downarrow r_2(3)&&\downarrow r_2(4)\\
&&&&&&\\
{\cal D}_{\R}^{0,0}(1)& \stackrel{D}{\lra} &({\cal D}^{0,1}\oplus {\cal D}^{1,0})_{\R}(1)&\stackrel{D}{\lra} &{\cal D}_{\R}^{1,1}(1)&
\stackrel{2 \overline \partial \partial}{\lra} &{\cal D}_{cl, \R}^{2,2}(2)
\end{array}
$$
where $D$ is the differential in ${\cal C}_{{\cal D}}(X, \R(2))$ 
$$
r_2(4): Y_2 \lms (2 \pi i)^2 \cdot \delta_{Y_2}; \qquad 
r_2(3): (Y_1,f) \lms 2 \pi i \cdot  \log|f|\delta_{Y_1}
$$
$$
r_2(2): f\wedge g \lms -\log|f|d i\arg g + \log|g|d i\arg f; \qquad 
r_2(1): \{f\}_2 \lms \widehat {\cal L}_2(f)
$$
To prove that we get a morphism of complexes we use Theorem \ref{CANHH}. The 
following 
argument is needed to check the commutativity of the second square.
The de Rham differential of the distribution $r_2(2)(f\wedge g)$ is 
$$
{\cal D}  \Bigl(-\log|f|d i\arg g + \log|g|d i\arg f\Bigr) = 
$$
$$
\pi_2(d \log f \wedge \log g) + 2 \pi i \cdot (\log|g| \delta(f) - \log|f| \delta(g)) 
$$ 
This  {\bf does not} coincide with $r_2(3) \circ \partial (f\wedge g)$, but the difference is
$$
({\cal D} \circ r_2(2) - r_2(3) \circ \partial) (f\wedge g) = \pi_2(d \log f \wedge \log g) \in \quad ({\cal D}^{0,2} \oplus {\cal D}^{2,0})_{\R}(1)
$$ 
Defining the differential $D$ on the second group of the complex 
${\cal C}_{{\cal D}}(X, \R(2))$ 
we take the de Rham differential  and throw away from it precisely these components. Therefore the middle square is commutative.

{\bf 5. The weight three}. The weight three motivic complex $\Gamma(X;3)$ is the total complex of the following bicomplex:
(the first group is in degree $1$)
$$
\begin{array}{ccccc}
{\cal B}_3(\Q(X))& \lra &{\cal B}_2(\Q(X)) \otimes \Q(X)^*& \lra &\Lambda^3\Q(X)^*\\
&&\downarrow &&\downarrow \\
&&\oplus_{Y_1 \in X^{(1)}}{\cal B}_2(\Q(Y_1))&\lra &\oplus_{Y_1 \in X^{(1)}}\Lambda^2\Q(Y_1)^*\\
&&&&\downarrow \\
&&&&\oplus_{Y_2 \in X^{(2)}}\Q(Y_2)^*\\
&&&&\downarrow \\
&&&&\oplus_{Y_3 \in X^{(3)}}\Q(Y_3)^*
\end{array}
$$
The Deligne complex 
${\cal C}_{{\cal D}}(X, \R(3))$ looks as follows:
$$
\begin{array}{ccccccc}
&&&&&&{\cal D}^{3,3}_{cl, \R}(3)\\
&&&&&2 \overline \partial \partial \nearrow &\\
{\cal D}^{0,2}  & \stackrel{\partial}{\lra} &{\cal D}^{1,2}  &\stackrel{\partial}{\lra} &{\cal D}^{2,2}  &&\\
\uparrow \overline \partial &&\uparrow \overline \partial &&\uparrow \overline \partial &&\\
{\cal D}^{0,1}  &\stackrel{\partial}{\lra} &{\cal D}^{1,1}  &\stackrel{\partial}{\lra} &{\cal D}^{1,2}  &&\\
\uparrow \overline \partial &&\uparrow \overline \partial &&\uparrow \overline \partial &&\\
{\cal D}^{0,0}  &\stackrel{\partial}{\lra} &{\cal D}^{1,0}  &\stackrel{\partial}{\lra} &{\cal D}^{2,0} &&
\end{array}
$$
We construct the regulator map $\Gamma(X;3) \lra {\cal C}_{{\cal D}}(X, \R(3))$
by setting 
$$
r_3(6): Y_3 \lms (2 \pi i)^3 \cdot \delta_{Y_3}; \qquad 
r_3(5): (Y_2,f) \lms (2 \pi i)^2 \cdot  \log|f|\delta_{Y_2}
$$
$$
r_3(4): (Y_1, f\wedge g) \lms 2 \pi i \cdot (-\log|f|d i\arg g + \log|g|d i\arg f )\delta_{Y_1}
$$
$$
r_3(3): (Y_1,\{f\}_2) \lms 2 \pi i \cdot\widehat {\cal L}_2(f)\delta_{Y_1}
$$
$$
r_3(3): f_1\wedge f_2 \wedge f_3 \lms {\rm Alt}_3\Bigl(\frac{1}{6}\log|f_1|d \log|f_2| \wedge d \log|f_3| + \frac{1}{2}\log|f_1|
d i\arg f_2 \wedge    d i\arg f_3\Bigr)
$$
$$
r_3(2): \{f\}_2 \otimes g \lms \widehat {\cal L}_2(f)d i \arg g - \frac{1}{3} \log|g| \cdot 
\Bigl(-\log |1-f| d\log |f| + \log |f| d\log |1-f|\Bigr)
$$
$$
r_3(1): \{f\}_3 \lms \widehat {\cal L}_3(f)
$$

{\bf 6. Classical polylogarithms and 
special values of $\zeta$-functions of algebraic varieties}. 
We conjecture that 
the polylogarithmic motivic complex $\Gamma(X;n)\otimes \Q$ should be 
quasiisomorphic to the weight $n$ motivic complex, and Beilinson's 
regulator map under this 
quasiisomorphism should be equal to the 
defined above regulator map on $\Gamma(X;n)$. 
This implies that the special values  
of $\zeta$-functions of algebraic varieties outside of the 
critical strip 
should be expressed via classical polylogarithms. 

The very special case of this conjecture when $X= {\rm Spec}(F)$ 
where $F$ is a number field is equivalent to Zagier's conjecture \cite{Z1}. 

The next interesting case 
is when $X$ is a curve over a number field.  
The conjecture in this case was elaborated in 
\cite{G3}, see also \cite{G8} for a survey. For example if $X$ is an elliptic curve 
this conjecture suggests that the special values 
$L(X,n)$ for $n \geq 2$ are expressed via certain generalized 
Eisenstein--Kronecker series. For $n=2$ this was 
discovered of S. Bloch \cite{Bl4}, and  
for $n=3$ it implies Deninger's conjecture \cite{De}.

A homomorphism
from $K$-theory to $H^i(\Gamma(X;n))\otimes \Q $ 
has been constructed in the following cases:

(1) $F$ is an arbitrary field, $n \leq 3$ (\cite{G1}, \cite{G2}, \cite{G5}]) 
and $n=4$, $i>1$ (to appear).

(2) $X$ is a curve over a number field, $n \leq 3$ (\cite{G5}) and $n=4$, $i>1$ 
(to appear). 

In all these cases we proved that this homomorphism 
followed by the regulator map 
on 
polylogarithmic complexes (when $F = \C(X)$ in (1)) coincides 
with Beilinson's regulator. 
This proves the difficult ``surjectivity'' property: 
the image of the regulator map on these polylogarithmic 
complexes {\it contains} 
the image of Beilinson's regulator map in the Deligne cohomology. 
The main ingredient of the proof is an explicit construction 
of the motivic cohomology class
$$
H^{2n}(BGL_{\bullet}, \Gamma(\ast; n))
$$
 
The results (2) combined with the results of R. de Jeu \cite{dJ1}-
\cite{dJ2} prove 
that the image of Beilinson's regulator map 
in the Deligne cohomology {\it coincides} with the image 
of the regulator map on polylogarithmic complexes in the case (2). 

{\bf 7. Coda: Polylogarithms on curves, \index{Polylogarithms on curves} Feynman integrals and special values 
of $L$-functions}. 
The classical polylogarithm functions admit the following generalization. 
Let $X$ be a regular complex algebraic curve. Let us assume first that 
$X$ is projective. Let us choose a metric on $X(\C)$. Denote by $G(x,y)$ the Green function 
provided by this metric. Then one can define 
real-valued functions $G_n(x,y)$ depending on a pair of points $x,y$ of $X(\C)$ with values in 
the complex vector space $S^{n-1}H_1(X(\C), \C)(1)$, see \cite{G10}, Chapter 9.1. 
By definition $G_1(x,y)$ is the Green function $G(x,y)$. It has a singularity at the diagonal, 
but provides 
a generalized function. 
For $n>1$ the function $G_n(x,y)$ is well-defined on $X(\C) \times X(\C)$. 
Let $X =E$ be an elliptic curve. Then the functions $G_n(x,y)$ are translation invariant, 
and thus reduced to a signle variable functions: $G_n(x,y) = G_n(x-y)$. 
The function $G_n(x)$ is given by the classical Kronecker-Eisenstein series. 
Finally, adjasting this construction to the case when $X = \C^*$, we get a 
function $G_n(x,y)$ invariant under the shifts on the group $\C^*$, i.e. one has 
$G_n(x,y) = G_n(x/y)$. The function $G_n(z)$  is given by 
the single-valued version (\ref{6.28.04.1}) of the classical $n$-logarithm function. 

Let us extend the function $G_n(x,y)$ by linearity to a function $G_n(D_1, D_2)$ 
depending on a pair of divisors on $X(\C)$. 
Although the function $G_n(x,y)$ depends on the choice of metric on $X(\C)$, 
the restriction of the function $G_n(D_1, D_2)$ to the subgroups of the 
degree zero divisors is independent of the choice of the metric.

Let $X$ be a curve over a number field $F$. 
In the Section 9.1 of \cite{G10} we proposed a conjecture which allows to express 
the special values $L({\rm Sym}^{n-1}H^1(X), n)$ via determinants whose entries are given by 
the values of the function $G_n(D_1, D_2)$, where $D_1, D_2$ are degree zero divisors on 
$(X(\overline F))$ invariant under the action of the Galois group ${\rm Gal}(\overline F/F)$.  
In the case when $X$ is an elliptic curve it boils down  
to the so-called elliptic analog of Zagier's conjecture, see \cite{Wil}, \cite{G11}. 
It has been proved for $n=2$ at \cite{GL}, and a part of this proof (minus surjectivity) 
can be transformed 
to the case of arbitrary $X$. 

It was conjectured in the Section 9.3 of \cite{G10} that the special values 
$L({\rm Sym}^{n-1}H^1(X), n+m-1)$, where $m \geq 1$, 
should be expressed similarly via the special values of the 
depth $m$ 
multiple polylogarithms on the curve $X$, defined in the Section 9.2 of loc. cit. 
One can show that in the case when $X=E$ is an elliptic curve, and $n = m = 2$, 
this conjecture reduces to Deninger's conjecture \cite{De} 
on $L(E,3)$, which has been proved in \cite{G3}  
An interesting  aspect of this story is that the multiple polylogarithms on curves  
\index{multiple polylogarithms on curves}  
are introduced via mathematically well defined Feynman integrals. 
This seems to be the first application of Feynman integrals in number theory. 

\section{Motivic Lie algebras and 
Grassmannian polylogarithms}

{\bf 1. Motivic Lie coalgebras and motivic complexes}. 
Beilinson conjectured that for 
an arbitrary field $F$ there exists an abelian $\Q$--category ${\cal M}_T(F)$ 
of mixed Tate motives over $F$. This category is supposed to be a mixed Tate category, see 
Section 8 of \cite{G12} for the background. Then the  Tannakian formalism 
implies that there exists a positively graded Lie coalgebra 
${\cal L}_{\bullet}(F)$ such that the category of finite 
dimensional graded comodules over ${\cal L}_{\bullet}(F)$ is naturally equivalent to 
the category ${\cal M}_T(F)$. 

Moreover ${\cal L}_{\bullet}(F)$ depends functorially on $F$. 
This combined with Beilinson's conjectural 
formula for the Ext groups in the category ${\cal M}_T(F)$ imply 
that the cohomology of this Lie coalgebra are computed 
by the formula
\begin{equation} \label{3.31.03.1}
H^i_{(n)}({\cal L}_{\bullet}(F)) = {\rm gr}^{\gamma}_{n}K_{2n-i}(F) \otimes \Q
\end{equation}
where $H^i_{(n)}$ means the degree $n$ part of $H^i$. 
This conjecture  provides a new point of view on 
algebraic $K$-theory of fields, suggesting and 
explaining several conjectures and results, see ch. 1 \cite{G2}.  

The degree $n$ part of the standard cochain complex 
$$
{\cal L}_{\bullet}(F) \stackrel{\delta}{\lra} \Lambda^2{\cal L}_{\bullet}(F) 
\stackrel{\delta \wedge {\rm Id} - {\rm Id}\wedge \delta}{\lra}
\Lambda^3{\cal L}_{\bullet}(F) \lra ...
$$
of the Lie coalgebra ${\cal L}_{\bullet}(F)$ is supposed to be quasiisomorphic to the 
weight $n$ motivic complex of $F$, providing yet another 
point of view on motivic complexes. 
For example the first four of the motivic complexes should look as follows:
$$
{\cal L}_1(F); \qquad {\cal L}_{2}(F) \lra \Lambda^2{\cal L}_1(F); \qquad {\cal L}_{3}(F) \lra 
{\cal L}_{2}(F) \otimes {\cal L}_1(F)\lra 
\Lambda^3{\cal L}_{1}(F)
$$
$$
{\cal L}_{4}(F) \lra {\cal L}_{3}(F)\otimes {\cal L}_1(F) \oplus \Lambda^2{\cal L}_2(F)\lra 
{\cal L}_{2}(F) \otimes \Lambda^2{\cal L}_1(F)\lra 
\Lambda^4{\cal L}_{1}(F)
$$
Comparing this with the formula (\ref{3.31.03.1}), 
the shape of complexes 
(\ref{3.31.03.2}) and (\ref{3.19.03.1}), and the known information relating 
their cohomology with algebraic $K$-theory we conclude the following.
One must have 
$$
{\cal L}_1(F) = F^*\otimes \Q; \quad {\cal L}_2(F) = B_2(F) \otimes \Q
$$
and we expect to have an isomorphism
$$
B_3(F)\otimes \Q \stackrel{\sim}{\lra} {\cal L}_3(F)
$$
Moreover the motivic complexes (\ref{3.31.03.2}) and (\ref{3.19.03.1}) are simply the degree $n$ parts of the standard cochain complex 
of the Lie coalgebra ${\cal L}_{\bullet}(F)$. 
More generally, the very existence for an arbitrary field $F$  of the 
motivic classical polylogarithms (\ref{cani2}) implies that 
one should have canonical homomorphisms
$$
l_n: {\cal B}_n(F) \lra {\cal L}_n(F)
$$
These homomorphisms are expected to satisfy the basic relation \ref{cani3}. Therefore 
the maps $l_n$ give rise to a canonical homomorphism from the weight $n$ polylogarithmic complex 
of $F$ to the degree $n$ part of the standard 
cochain complex of ${\cal L}_{\bullet}(F)$. 
\begin{equation}\label{6.30.04.10}
\begin{array}{ccccc}
{\cal B}_n(F)& \lra &{\cal B}_{n-1}(F)\otimes F^*& \lra ... \lra & \Lambda^nF^*\\
&&&&\\
\downarrow l_n&&\downarrow l_{n-1}\wedge l_1&&\downarrow =\\
&&&&\\
{\cal L}_n(F)&\lra &\Lambda^2_{(n)}{\cal L}_{\bullet}(F)&\lra ... \lra 
&\Lambda^n{\cal L}_1(F)\\
\end{array}
\end{equation}
where $\Lambda^2_{(n)}$ denotes the degree $n$ part of $\Lambda^2_{(n)}$. 
 For $n=1,2,3$, these maps, combined with the ones 
(\ref{4.15.03.1}), lead to the 
maps above. For $n>3$ the  map of complexes (\ref{6.30.04.10}) will not be an isomorphism. 
The conjecture that it is a quisiisomorphism is equivalent to the 
Freeness conjecture for the Lie coalgebra ${\cal L}_{\bullet}(F)$, see 
\cite{G1}-\cite{G2}. 

Therefore we have two different points of view on the groups $B_n(F)$ 
for $n =1,2,3$: according to one of them they are 
the particular cases of the groups ${\cal B}_n(F)$, and 
according to the other they provide an explicit 
computation of the first three of the groups ${\cal L}_n(F)$. 
It would be very interesting to find an explicit  construction of the groups 
${\cal L}_n(F)$ for $n>3$ generalizing the definition of the groups $B_n(F)$. 
More specifically, 
we would like to have a ``finite  dimensional'' construction 
of all vector spaces ${\cal L}_{n}(F)$, i.e. for every $n$ 
there should exist a finite number of finite 
dimensional algebraic  varieties $X_n^i$ and $R_n^j$ such that 
$$
{\cal L}_n(F) = {\rm Coker}\left(\oplus_j \Q[R_n^j(F)] \lra \oplus_i \Q[X_n^i(F)]\right)
$$
Such a construction would be provided by 
the scissor congruence 
motivic Hopf algebra of $F$ (\cite{BMS}, \cite{BGSV}). However so far 
the problem in the definition of the coproduct 
in loc. cit. for non generic generators has not been resolved. 
A beautiful construction of the motivic Lie coalgebra of a field $F$ 
was suggested by Bloch and Kriz in \cite{BKr}. However 
it is not finite dimensional in the above sense. 

{\bf 2. Grassmannian approach of the motivic Lie coalgebra}. 
We suggested in \cite{G5} that there should exist a 
construction of the Lie coalgebra ${\cal L}_{\bullet}(F)$ such that 
the variety $X_n$ is the variety of configurations of 
$2n$ points in $P^{n-1}$ and the relations varieties 
$R_n^j$ are provided by the functional equations for the 
{\it motivic Grassmannian $n$-logarithm}. 
Let us explain this in more detail.

Let $\widetilde {\cal L}_{n}(F)$ be the free abelian group generated by 2n-tuples of
points $(l_1,...,l_{2n})$ in 
 ${{ P}}^{n-1}(F)$ subject to the following relations:

1)$(l_1,...,l_{2n}) = 
(gl_1,...,gl_{2n})$ for any $g \in PGL_{n}(F)$.

2) $(l_1,...,l_{2n} ) = (-1)^{|\sigma|} (l_{\sigma
(1)},...,l_{\sigma (2n)})$
for any  permutation $\sigma \in S_{2n}$.

3) for any $2n+1$ points $(l_0,...,l_{2n})$ in $P^{n-1}(F)$ one has
$$ 
\sum_{i=0}^{2n}(-1)^{i}(l_0,...\hat
l_i,...,,l_{2n}) =0
$$ 

4) for any $2n+1$ points $(l_0,...,l_{2n})$ in ${{P}}^{n}(F)$ one has
$$ 
\sum_{i=0}^{2n}(-1)^{i}(l_i|l_0,...\hat
l_i,...,,l_{2n}) =0
$$ 

We conjecture that ${\cal L}_{n}(F)$ is a quotient of 
$\widetilde {\cal L}_{n}(F)$. It is a nontrivial quotient already for $n=3$. 
Then to define the Lie coalgebra ${\cal L}_{\bullet}(F)$ one needs 
to produce a cobracket
$$
\delta: {\cal L}_{n}(F) \lra\oplus_{i} {\cal L}_{i}(F)\wedge {\cal L}_{n-i}(F)
$$
Here is how to do this in the first nontrivial case, $n=4$.

Let us define   a homomorphism 
$$
\widetilde {\cal L}_4(F) 
\quad \stackrel{ \delta }{\longrightarrow} \quad B_3(F) \otimes F^{\ast}
\quad \oplus  \quad B_2(F)\wedge  B_2(F)   
$$
by setting  $ \delta = ( \delta_{3,1},    \delta_{2,2})$ where
$$
 \delta_{3,1}(l_1,...,l_8) : = 
-  \frac{1}{9}{\rm Alt}_8\Bigl( \Bigl( 
r_3(l_1|,l_2,l_3,l_4;l_5,l_6,l_7)  +   
\{r(l_1 l_2| l_3,  l_6, l_4,  l_5)\}_3 -
$$
$$
- 
\{r(l_1 l_2| l_3,  l_5, l_4,  l_6)\}_3 
\Bigr) \wedge \Delta (l_5,l_6,l_7,l_8)\Bigr)
\in B_3(F) \wedge F^*
$$
$$
 \delta_{2,2}(l_1,...,l_8) : =  
$$
$$
\frac{1}{7} \cdot {\rm Alt}_8
\Bigl( \{r_2(l_1,l_2|l_3,l_4,l_5,l_6)\}_2\wedge 
\{r_2(l_3, l_4|l_1,l_2,l_5,l_7)\}_2\Bigr) \in \Lambda^2 B_2(F)
$$
These formulae are obtained by combining the definitions 
at p. 136-137 and p 156 of \cite{G5}. The following key result is 
Theorem 5.1 in loc. cit. 

\begin{theorem} \label{mmmar} a) The homomorphism $ \delta$ kills the relations 1) - 4).

b) The following composition equals to  zero:
$$
\widetilde {\cal L}_4(F) \stackrel{ \delta }{\longrightarrow} B_3(F) \otimes F^{\ast}
\quad \oplus \quad  B_2(F)\wedge  B_2(F)  \stackrel{\delta}{\lra} B_2(F) \otimes \Lambda^2F^{\ast} 
$$
\end{theorem}
Here the second differential is defined  by the Leibniz rule, 
using the differentials in complexes (\ref{3.31.03.2}) and (\ref{3.19.03.1}). 

Taking the "connected component" of zero in ${\rm Ker} \delta $ 
(like we did in  s. 1 ch. 4 above, or in \cite{G1})
 we should get the set of defining relations for the group ${\cal L}_4(F)$. 
An explicit construction of them is not known.

{\bf 4. The motivic Grassmannian tetralogarithm}. 
Let $F$ be a function field on a complex variety $X$. 
There is a morphism of complexes
$$
\begin{array}{ccccc}
 B_3(F) \otimes F^{\ast}
\quad \oplus \quad  B_2(F)\wedge  B_2(F)& \stackrel{\delta }{\longrightarrow} &B_2(F) \otimes \Lambda^2F^*& \stackrel{\delta }{\longrightarrow}& \Lambda^4F^*\\
&&&&\\
\downarrow R_4(2)&&\downarrow R_4(3)&&\downarrow R_4(4)\\
&&&&\\
{\cal A}^1(Spec F)& \stackrel{d}{\longrightarrow}& {\cal A}^2(Spec F)&\stackrel{d}{\longrightarrow}& {\cal A}^3(Spec F)
\end{array}
$$
extending the homomorphism $r_4(*)$ from Chapter 5. 
Namely, $R_4(*) = r_4(*)$ for $* = 3,4$ and $R_4(2) = (r_4(2), r'_4(2))$ where 
$$
r'_4(2): \Lambda^2B_2(F) \lra {\cal A}^1(Spec F)
$$
$$
\{ f \}_{2} \wedge \{ g \}_{2}  \mapsto   \frac{1}{3}\cdot \Bigl(   \widehat{{\cal L}}_{2} (g )   \cdot \alpha(1-f,f) - \widehat{{\cal L}}_{2} (f )   \cdot \alpha(1-g,g) \Bigr)
$$
  It follows from Theorem \ref{mmmar} that the composition 
$R_4(2) \circ \delta(l_1, ..., l_8)$ is a closed 1-form on 
the configuration space of $8$ points in $\CP^3$. 
One can show (Proposition 5.3 in \cite{G5}) that 
integrating this $1$-form we get a single valued function 
on the configuration space, denoted ${\cal L}_4^M$ and called 
the {\it motivic Grassmannian tetralogarithm}. 
It would be very interesting to compute 
the difference ${\cal L}_4^M - {\cal L}^{\cal L}_4$, similarly to the
formula (\ref{3.31.03.5}) in the case $n=3$.   We expect that 
is expressed as a sum of products of functions 
${\cal L}_3$, ${\cal L}_2$ and $\log|*|$.

\begin{theorem} \label{5.2th} a) There exists a canonical  map 
$$
K^{[3]}_7(F)\otimes {\Q} \quad \lra \quad {\rm Ker} \Bigl(\widetilde {\cal L}_4(F)  \stackrel{\delta}{\lra} B_3(F) \otimes F^* \oplus \Lambda^2B_2(F)\Bigr)_{\Q}
$$

b) In the case $F = {\C}$ 
the composition 
$
K_7({\C} ) \lra \widetilde {\cal L}_4({\C})  \stackrel{{\cal L}^M_4}{\lra} \R
$ 
coincides with a nonzero rational multiple of the Borel regulator map. 
\end{theorem}

The generalization of the above picture to the case $n>4$ is unknown. 
It would be very interesting at least to define 
the motivic Grassmannian polylogarithms via
 the Grassmannian polylogarithms. 

For $n=3$ the motivic Grassmannian trilogarithm ${\cal L}^M_3$ 
is given by the first term of 
the formula (\ref{3.31.03.5}). It is known (\cite{GZ}) that 
${\cal L}^{\cal L}_3$ does not provide a homomorphism 
${\cal L}_3(\C) \lra \R$ since the second term in 
(\ref{3.31.03.5}) does not have this property. 
Indeed, the second term in (\ref{3.31.03.5}) vanishes 
on the special configuration of $6$ points in $P^2$, but does not vanish 
at the generic configuration. On the other hand 
the defining relations in ${\cal L}_3(F)$ allow to express any 
configuration as a linear combination of the special ones. 
 
We expect the same situation in general: for $n>2$
the Grassmannian $n$--logarithms should not satisfy all 
the functional equations for the motivic Grassmannian 
$n$--logarithms. It would be very interesting to find a 
conceptual explanation of this surprising phenomena.


\begin{thebibliography}{BGSV}

\bibitem[B1]{B1} Beilinson A.A: {\it Higher regulators and values of 
$L$-functions}, VINITI, 24 (1984), 181--238. 

\bibitem[B2]{B2} Beilinson A.A: {\it Height pairing between 
algebraic cycles},  Springer Lect. Notes in Math. 1289, 1-25. 

\bibitem[B3]{B3} Beilinson A.A: {\it Notes on absolute Hodge cohomology}. 
Applications of algebraic $K$-theory to algebraic geometry and number theory, 
Part I, II (Boulder, Colo., 1983),
   35--68, Contemp. Math., 55, Amer. Math. Soc., Providence, RI, 1986. 


\bibitem[BD1]{BD1} Beilinson A.A, Deligne P.: {\it Interprétation motivique de la conjecture de Zagier reliant polylogarithmes et régulateurs}, Motives 
(Seattle, WA, 1991), 97--121, Proc. Sympos. Pure Math., 55, Part 2, 
Amer. Math. Soc., Providence, RI, 1994.

\bibitem[BD2]{BD2} Beilinson A.A, Deligne P.: 
{\it Polylogarithms}, Manuscript, 1992. 

\bibitem[BMS]{BMS} Beilinson A.A., MacPherson R., Schechtman, 
V.V: {\it Notes on motivic cohomology}, Duke Math.\ J.\ 54 (1987), 
679--710. 

\bibitem[BGSV]{BGSV} Beilinson A.A., Goncharov A.A., Schechtman V.V., 
Varchenko A.N.: {\it Aomoto dilogarithms, mixed Hodge structures and 
motivic cohomology}, 
The Grothendieck Festschrift, Birkhauser, vol 86, 1990, p. 135-171.



\bibitem[Bl1]{Bl1} Bloch S: {\it Algebraic cycles and higher K-theory}, Adv. in
Math., (1986), v. 61, 267--304.

\bibitem[Bl2]{Bl2} Bloch S: {\it The moving lemma for higher Chow groups},
 J. Alg. geometry. 3 (1994), 537-568.

\bibitem[Bl3]{Bl3} Bloch S: {\it Algebraic cycles and the 
Beilinson conjectures}, Cont. Math. 58 (1986) no 1, 65-79.

\bibitem[Bl4]{Bl4} Bloch S: {\it Higher regulators, algebraic K-theory and zeta functions of elliptic curves}, Irvines lectures, CRM Monograph series, vol. 11, 
2000. 

\bibitem[Bl5]{Bl5} Bloch, S: {\it Applications of the dilogarithm function in algebraic $K$-theory and algebraic geometry},   
Proceedings of the International Symposium on Algebraic Geometry (Kyoto Univ.,
   Kyoto, 1977), pp. 103--114.


\bibitem[BK]{BK} Bloch S, Kato K.: {\it $L$-functions and 
Tamagawa numbers of motives}, In: 
The Grothendieck Festschrift, Birkha\"user, vol 86, 1990.

\bibitem[BKr]{BKr}Bloch S, Kriz I.: {\it Mixed Tate motives},
Ann. Math. 140 (1994), 557-605. 

\bibitem[Bo1]{Bo1} Borel A: {\it Stable real cohomology of arithmetic groups}, Ann. Sci. Ec. Norm. Super., (4) 7 (1974), 235--272.
\
\bibitem[Bo2]{Bo2} Borel A: {\it Cohomologie de $SL_{n}$ et valeurs de 
fonctions z\^{e}ta aux points entiers},  Ann. Sc. Norm. 
Sup. Pisa 4 (1977), 613--636.


\bibitem[Bu]{Bu} Burgos, J. I.: {\it 
Arithmetic Chow rings and Deligne-Beilinson cohomology},
 J. Algebraic Geom. 6 (1997), no. 2, 335--377. 

\bibitem[D]{D} Deligne P.: {\it 
Valeurs de fonctions $L$ et periodes d'integrales},
In: AMS Proc. Symp. Pure Math., vol. 33, part 2, (1979). 313-346. 

\bibitem[D2]{D2} Deligne P.: {\it Le groupe fundamental de la droit projective
    moins trois points}, In:  Galois groups over $\Q$, MSRI Publ. 16, 79-313,
    Springer Verlag, 1989. 
. 

\bibitem[DG]{DG} Deligne P., Goncharov A.: {\it Groupes fondamentaux motiviques de Tate mixte}, Ann. Sci. Ecole Norm. Sup., To appear. arXive: math.NT/0302267.

\bibitem[De]{De} Deninger Ch. 
{\it Higher order operatons in Deligne cohomology}, Inventiones Math.
122, N1 (1995) 289-316.

\bibitem[Du]{Du} Dupont J.: {\it Simplicial De Rham cohomolgy and
characteristic classes of flat bundles}, Topology 15, 1976, p. 233-245.

\bibitem[GGL]{GGL} Gabrielov A., Gelfand I.M., Losik M.: 
{\it Combinatorial computation of characteristic classes}, I, II,
Funct. Anal. i ego pril. 9 (1975) 12-28, ibid 9 (1975) 5-26.


\bibitem[GaZ]{GaZ} Gangl H., Zagier D.: 
{\it  Classical and elliptic polylogarithms and special values of L-series},
The arithmetic and geometry of algebraic cycles (Banff, AB, 1998), 561--615, NATO
   Sci. Ser. C Math. Phys. Sci., 548, Kluwer Acad. Publ., Dordrecht, 2000. 

\bibitem[GM]{GM} Gelfand I.M., MacPherson R: {\it Geometry in 
Grassmannians and a generalisation of the dilogarithm}, Adv. 
in Math., 44 (1982) 279--312. 

 \bibitem[GS]{GS} Gillet, H.; Soul\'e, Ch.:  
{\it Arithmetic intersection theory}, 
Inst. Hautes ŽÉtudes Sci. Publ. Math. No. 72, (1990), 93--174 (1991). 

\bibitem[G1]{G1} Goncharov A.B: {\it Geometry of configurations, 
polylogarithms and motivic cohomology},  
Adv. in Math., 114, (1995), 197--318. 

\bibitem[G2]{G2} Goncharov A.B.: {\it Polylogarithms and motivic Galois
groups}, in Symp. in Pure Math., v. 55, part 2, 1994, 43 - 97.



\bibitem[G3]{G3}
Goncharov, A.B.: {\it Deninger's conjecture on special values of
$L$-functions of elliptic curves at $s=3$},  J. Math. Sci. 81 (1996), N3,  2631-2656,  
alg-geom/9512016. 


\bibitem[G4]{G4} Goncharov A.B.: {\it  Chow polylogarithms and
regulators}, Math. Res. Letters, 2, (1995), 99-114.



\bibitem[G5]{G5}
Goncharov, A.B.: {\it Geometry of trilogarithm and the motivic Lie 
algebra of a field}, Proceedings of the conference ``Regulators in 
arithmetic and geometry'', Progress in Mathematics, vol 171, Birkh\"auser 
(2000), 127-165. math.AG/0011168.

\bibitem[G6]{G6}
Goncharov, A.B.: {\it Explicit regulator maps on the polylogarithmic 
motivic complexes}, In:
`` Motives, polylogarithms and Hodge theory''. Part I,
 F. Bogomolov, L.  Katzarkov editors. 
International Press. 245-276. 
 math.AG/0003086.


\bibitem[G7]{G7}
Goncharov, A.B.: {\it Polylogarithms, 
regulators, and Arakelov motivic complexes}, To appear in JAMS, 2005,   
math.AG/0207036. 

\bibitem[G8]{G8} 
Goncharov, A.B.: {\it Polylogarithms in arithmetic and geometry},  
Proceedings of the International Congress of
  Mathematicians, Vol. 1, 2 (Zurich, 1994), 374--387, 
Birkhauser, Basel, 1995.
 
\bibitem[G9]{G9} 
Goncharov, A.B.: {\it Periods and mixed motives},  Preprint 
math.AG/0202154.  

\bibitem[G10]{G10} 
Goncharov, A.B.: {\it Multiple zeta-values, Galois groups, 
and geometry of modular varieties}, European Congress of Mathematics, Vol. I (Barcelona, 2000), 361--392, Progr. Math., 201,
   Birkhauser, Basel, 2001. math.AG/0005069.
 
\bibitem[G11]{G11} Goncharov, A.B.: {\it Mixed elliptic motives}. Galois 
representations in arithmetic algebraic geometry (Durham, 1996), 147--221, 
London Math. Soc. Lecture Note Ser., 254, Cambridge Univ. Press,
   Cambridge, 1998. 

\bibitem[G12]{G12} Goncharov, A.B.: {\it Galois symmetries of fundamental groupoids and noncommutative geometry}, To appear in Duke Math. Journal, 2005,  math.AG/0208144. 

\bibitem[GL]{GL} 
Goncharov, A.B., Levin A.M.: {\it Zagier's conjecture on $L(E,2)$}. 
Invent. Math. 132 (1998), no. 2, 393--432. 

\bibitem[GoM]{GoM} {\it 
Multiple $\zeta$-motives and moduli spaces $\overline{\cal M}\sb {0,n}$.}  Compos. Math.  140  (2004),  no. 1, 1--14.

\bibitem[GZ]{GZ} 
Goncharov A. B., Zhao J.: {\it The Grassmannian trilogarithms}, 
Compositio Math. 127 (2001), no. 1, 83--108. math.AG/0011165.

 

\bibitem[HM]{HM} Hain R, MacPherson R: {\it Higher Logarithms}, 
Ill. J.\ of Math,, vol. 34, (1990) N2,  392--475. 

\bibitem[HW]{HW} Huber A., Wildeshaus J. {\it Classical motivic polylogarithms 
according to Beilinson and Deligne},  Doc. Math. 3 (1998), 27--133 (electronic). 

\bibitem[dJ1]{dJ1} de Jeu, Rob: {\it On $K^{(3)}_4$ of 
curves over number fields}, 
Invent. Math. 125 (1996), N3, 523-556. 

\bibitem[dJ2]{dJ2} de Jeu, Rob: {\it Towards the regulator 
formulas for curves over number fields}, Compositio Math. 124
  (2000), no. 2, 137--194. 

\bibitem[dJ3]{dJ3} de Jeu, Rob: {\it Zagier's conjecture 
and wedge complexes 
in algebraic K-theory}, Comp. Math., 96, N2, (1995), 197-247.

\bibitem[Ha]{Ha} Hanamura M.: {\it Mixed motives and algebraic cycles}, I. 
Math. Res. Lett. 2 (1995), no. 6, 811--821. 

\bibitem[HM1]{HM1} Hanamura M, MacPherson R: {\it Geometric construction
of polylogarithms}, Duke Math. J. 70 ( 1993)481-516. 

\bibitem[HM2]{HM2} Hanamura M,  MacPherson R.: 
{\it Geometric construction of polylogarithms, II}, 
Progress in Math. vol. 132 (1996) 215-282. 

\bibitem[Le]{Le} Levin A.M.: {\it Notes on $\R$-Hodge-Tate sheaves}, 
Preprint MPI 2001. 

\bibitem[Lev]{Lev} Levine, M: {\it Mixed motives}, 
Mathematical Surveys and Monographs, 57. American Mathematical Society, Providence, RI, 1998

\bibitem[Li1]{Li1} Lichtenbaum S: {\it Values of zeta functions, 
\'etale cohomology, and algebraic K-theory}, in Algebraic K-theory II, 
Springer LNM, 342, 1973, 489-501. 

\bibitem[Li2]{Li2} Lichtenbaum S: {\it Values of zeta functions at  
non-negative integers}, Lect. Notes in Math., 1068, Springer Verlag,
1984, 127-138. 

\bibitem[MW]{MW} Mazur B, Wiles A.: {\it Class fields of abelian 
extensions of $\Q$}, Invent. Math. 76 (1984), no. 2, 179--330.

\bibitem[Mi]{Mi} Milnor J. Introduction to algebraic K-theory. 
Princeton Univ. Press. 1971. 

\bibitem[N]{N} Nekov\'a\v{r} J.: {\it Beilinson's conjectures}. Motives (Seattle, WA, 1991), 537--570, Proc. Sympos. Pure Math., 55, Part 1, Amer. Math. Soc., Providence, RI, 1994.
   
\bibitem[RSS]{RSS}  Rapoport M., Schappacher N.  Schneider P.:
 Beilinson's conjectures on special values of $L$-functions. Perspectives in Mathematics, 4.
   Academic Press, Inc., Boston, MA, 1988.

\bibitem[R]{R}  Ramakrishnan D.: {\it Regulators, algebraic cycles, and values of $L$-functions}. Algebraic $K$-theory and algebraic number theory (Honolulu, HI, 1987), 183--310,
   Contemp. Math., 83, Amer. Math. Soc., Providence, RI, 1989.

                                              
\bibitem[S]{S} Soule C.,  D. Abramovich, J. F. Burnol, J. K. Kramer: 
Lectures on Arakelov geometry.  Cambridge University Press, 
1992. 

\bibitem[Sch]{Sch} Scholl A.: {\it Integral elements in $K$-theory and products of modular curves}. 
The arithmetic and geometry of algebraic cycles (Banff, AB, 1998), 467--489, 
NATO Sci. Ser. C Math. Phys. Sci., 548. 

\bibitem[Su]{Su} Suslin A.A.: {\it $K_3$ of the field and Bloch's group}, 
Proc. of the Steklov Institute of Math., 1991, issue 4, 217-240. 

\bibitem[T]{T} Tate J.: {\it Symbols in arithmetics} 
Actes du Congr\'es Int. Math. vol. 1, Paris, 1971, 201-211. 

\bibitem[V]{V} Voevodsky V.: 
{\it Triangulated category of motives over a field}, in 
Cycles, transfers, and motivic homology theories. 
Annals of Mathematics Studies,
143. Princeton University Press, Princeton, NJ, 2000. 

\bibitem[Wa]{Wa} Wang Q.: {Construction of the motivic multiple polylogarithms}, Ph. D. thesis, Brown University, April 2004. 


\bibitem[W]{W} Wiles A.: {\it The Iwasawa conjecture for totally 
real fields}, Ann. of Math. (2) 131 (1990), no. 3, 493--540. 

\bibitem[Wil]{Wil} Wildeshaus J.: {\it On an elliptic analogue of Zagier's conjecture},
Duke Math. J. 87 (1997), no. 2, 355--407.

\bibitem[Z1]{Z1} Zagier D: {\it Polylogarithms, Dedekind zeta 
functions and the algebraic $K$-theory of fields}, Arithmetic         
 Algebraic Geometry (G.v.d.Geer, F.Oort, J.Steenbrink, eds.),        
  Prog. Math., Vol 89, Birkhauser, Boston, 1991, pp. 391--430. 

\end{thebibliography}
\end{document}